\documentclass[11pt,letterpaper]{article}

\usepackage{amsfonts, amsmath, amssymb, amscd, amsthm, graphicx}

\newtheorem{thm}{Theorem}[section]
\newtheorem{cor}[thm]{Corollary}
\newtheorem{lem}[thm]{Lemma}
\newtheorem{prop}[thm]{Proposition}

\theoremstyle{definition}
\newtheorem{defn}[thm]{Definition}
\theoremstyle{remark}

\newfont{\eufm}{eufm10}

\newcommand{\G }{\Gamma (G, X\cup \mathcal H)}
\newcommand{\Ga }{\Gamma (G, \mathcal A)}
\newcommand{\dxh }{dist_{X\cup\mathcal H}}

\newcommand{\Hl }{\{ H_\lambda \} _{\lambda \in \Lambda } }

\newcommand{\e }{\varepsilon }
\newcommand{\N }{\mbox{\eufm N}}

\renewcommand{\P }{\mathcal P}
\renewcommand{\phi }{{\rm\bf Lab\, }}
\renewcommand{\kappa }{\varkappa}
\renewcommand{\AA }{Area^{rel}_{G(\N )}}
\begin{document}

\title{Peripheral fillings of relatively hyperbolic groups}
\author{D. V. Osin \thanks{This work has been partially supported by the NSF grant DMS-0605093 and by the
Russian Fund for Basic Research grant $\# $05-01-00895. }}
\date{}%

\maketitle

\begin{abstract}
In this paper a group theoretic version of Dehn surgery is studied.
Starting with an arbitrary relatively hyperbolic group $G$ we define
a {\it peripheral filling procedure}, which produces quotients of
$G$ by imitating the effect of the Dehn filling of a complete finite
volume hyperbolic 3--manifold $M$ on the fundamental group $\pi
_1(M)$. The main result of the paper is an algebraic counterpart of
Thurston's hyperbolic Dehn surgery theorem. We also show that
peripheral subgroups of $G$ `almost' have the Congruence Extension
Property and the group $G$ is approximated (in an algebraic sense)
by its quotients obtained by peripheral fillings.

\medskip


\smallskip

\end{abstract}

\section{Introduction}

Let $M$ be a compact orientable 3--manifold with finitely many
toric boundary components $T_1, \ldots , T_k$. Topologically
distinct ways to attach a solid torus to $T_i$ are parameterized
by {\it slopes} on $T_i$, i.e., isotopy classes of unoriented
essential simple closed curves in $T_i$. For a collection $\sigma
=(\sigma _1, \ldots , \sigma _k)$, where $\sigma _i$ is a slope on
$T_i$, the {\it Dehn filling } $M(\sigma )$ of $M$ is the manifold
obtained from $M$ by attaching a solid torus $\mathbb D^2\times
\mathbb S^1$ to each boundary component $T_i$ so that the meridian
$\partial \mathbb D^2$ goes to a simple closed curve of the slope
$\sigma _i$. The fundamental theorem of Thurston \cite{Th} asserts
that if $M-\partial M$ admits a complete finite volume hyperbolic
structure, then the resulting closed manifold $M(\sigma )$ is
hyperbolic provided $\sigma $ does not contain slopes from a fixed
finite set.

Given a subset $S$ of a group $G$, we denote by $\langle S\rangle
^G$ the normal closure of $S$ in $G$. Clearly, $$\pi _1(M(\sigma
))=\pi _1(M)/\langle x_1, \ldots , x_k\rangle ^{\pi _1(M)},$$ where
$x_i\in \pi_1 (T_i)\le \pi _1(M)$ is the element corresponding to
the slope $\sigma _i$. Thus Thurston's theorem implies the following
group theoretic result:

{\it Let $G$ be the fundamental group of a complete finite volume
hyperbolic $3$--manifold, $H_1, \ldots , H_k$ the cusp subgroups of
$G$. Then there exists a finite subset $\mathcal F$ of $G$ such that
for any collection of (primitive) elements $x_i\in H_i\setminus
\mathcal F$, the quotient group $G/\langle x_1, \ldots , x_k\rangle
^G$ is (word) hyperbolic.}

In our paper we generalize this result in two directions. First
instead of the class of fundamental groups of complete finite volume
hyperbolic manifolds we consider its far--reaching generalization,
the class of {\it relatively hyperbolic groups}. Secondary instead
of single elements $x_i\in H_i$ we deal with normal subgroups
generated by arbitrary subsets of the cusp subgroups.

Recall that the notion of relative hyperbolicity was introduced in
group theory by Gromov in \cite{Gro} and since then it has been
elaborated from different points of view \cite{Bow,DS,F,RHG}. Here
we mention some examples and refer the reader to the next section
for the precise definition of relatively hyperbolic groups.

\begin{itemize}

\item If $M$ is a complete finite--volume manifold of pinched
negative sectional curvature, then $\pi _1(M)$ is hyperbolic with
respect to the collection of the cusp subgroups \cite{Bow,F}.

\item Any (word) hyperbolic group $G$ is hyperbolic relative to
the trivial subgroup.

\item Geometrically finite convergence groups acting on non--empty
perfect compact metric spaces are hyperbolic relative to the set
of the maximal parabolic subgroups \cite{Y}.

\item Free products of groups and their small cancellation
quotients, as defined in \cite{LS}, are hyperbolic relative to the
factors \cite{RHG}.

\item Fundamental groups of finite graphs of groups with finite
edge groups are hyperbolic relative to the vertex groups \cite{Bow}.
In particular, according to the famous Stallings Theorem
\cite{Stall}, any group with infinite number of ends carries such a
relatively hyperbolic structure.

\item Finitely generated groups acting freely on $\mathbb R^n$--trees are
hyperbolic relative to the maximal non--cyclic abelian subgroups
\cite{Gui}. This class of examples includes limit groups studied by
Kharlampovich, Myasnikov \cite{KM}, and independently by Sela
\cite{Sel} in their solutions of the famous Tarskii problem.
\end{itemize}

In the context of relatively hyperbolic groups, the algebraic
analogue of Dehn filling is defined as follows. Suppose that $\Hl $
is a collection of subgroups of a group $G$. To each collection $\N
=\{ N_\lambda \} _{\lambda \in \Lambda }$, where $N_\lambda $ is a
normal subgroup of $H_\lambda $, we associate the quotient group
\begin{equation} \label{GN}
G(\N ) = G/\left\langle \mbox{$\bigcup_{\lambda \in \Lambda }$}
N_\lambda \right\rangle ^G .
\end{equation}

Our main result is the following.

\begin{thm}\label{CEP}
Suppose that a group $G$ is hyperbolic relative to a collection of
subgroups $\Hl $. Then there exists a finite subset $\mathcal F$ of
non--trivial elements of $G$ with the following property. Let $\N
=\{ N_\lambda \}_{\lambda \in \Lambda }$ be a collection of
subgroups $N_\lambda\lhd H_\lambda $ such that $N_\lambda \cap
\mathcal F=\emptyset $ for all $\lambda \in \Lambda $. Then:
\begin{enumerate}

\item[1)]  For each $\lambda\in \Lambda $, the natural map $H_\lambda
/N_\lambda \to G(\N )$ is injective.

\item[2)] The quotient group $G(\N )$ is hyperbolic relative to the
collection $\{H_\lambda /N_\lambda \} _{\lambda \in \Lambda }$.

\end{enumerate}
Moreover, for any finite subset $S\subseteq G$, there exists a
finite subset $\mathcal F(S)$ of non--trivial elements of $G$ such
that the restriction of the natural homomorphism $G\to G(\N )$ to
$S$ is injective whenever  $N_\lambda \cap \mathcal F(S)=\emptyset $
for all $\lambda \in \Lambda $.
\end{thm}

It is worthwhile to notice that the theorem applies to general (not
necessarily finitely generated) relatively hyperbolic groups. In
case the group $G$ is finitely generated, the condition $N_\lambda
\cap \mathcal F=\emptyset $ simply means that the subgroups
$N_\lambda $ contain no non--trivial elements of small (word)
length.

Our proof is purely combinatorial and extensively uses techniques
related to van Kampen diagrams over group presentations. Many ideas
used in the proof go back to methods developed by Alexander
Olshanskii in his geometric solution of the Burnside problem
\cite{Ols-book,Ols92}. After this paper was submitted, another proof
of Theorem \ref{CEP} in the particular case when the group $G$ is
torsion free and finitely generated was published by Daniel Groves
and Jason Manning in arXiv \cite{MG}. Later in \cite{MG1} they
showed that their method works for infinitely generated groups as
well.

Recall that if a finitely generated group $G$ is hyperbolic relative
to a collection of hyperbolic subgroups, then $G$ is a hyperbolic
group itself \cite{F,RHG}. The following corollary may be considered
as a generalization of the group theoretic version of Thurston's
hyperbolic Dehn surgery theorem. Indeed in case $G$ is a fundamental
group of a complete finite volume hyperbolic $3$--manifold, all cusp
subgroups are isomorphic to $\mathbb Z\oplus\mathbb Z$ and for any
non--trivial element $x\in \mathbb Z\oplus\mathbb Z$, the quotient
$\mathbb Z\oplus\mathbb Z/\langle x\rangle $ is hyperbolic.

\begin{cor}\label{hypquot}
Under the assumptions of Theorem \ref{CEP}, suppose in addition that
$G$ is finitely generated and $H_\lambda /N_\lambda $ is hyperbolic
for each $\lambda\in \Lambda $. Then $G(\N )$ is hyperbolic.
\end{cor}

On the other hand, Theorem \ref{CEP} can also be applied to
manifolds of higher dimension. Indeed let $M$ be a complete finite
volume hyperbolic $n$--manifold with cusp ends $E_1, \ldots , E_k$.
For simplicity we assume that each cusp end is homeomorphic to
$T^{n-1}\times \mathbb R^{+}$, where $T^{n-1}$ is an
$(n-1)$-dimensional torus. Let $1\le l\le k$. For each $E_i$, $1\le
i\le l$, we fix a torus $T_i^{n-1}\subset E_i$ and a closed simple
curve $\sigma _i$ in $T_i^{n-1}$. We now perform Dehn filling on the
collection of cusps $E_i$, $1\le i\le l$, by attaching a solid torus
$\mathbb D^2\times T^{n-2}$ onto $T^{n-1}_i$ via a homeomorphism
sending $\mathbb S^1=\partial \mathbb D^2$ to $\sigma _i$. Let
$\sigma =(\sigma _1, \ldots , \sigma _l)$. The topological type of
the resulting manifold $M(\sigma )$ depends only on the homotopy
class of unoriented curves  $\sigma _i$, i.e., on $[\pm \sigma
_i]\in \pi_1(T^{n-1}_i)$ \cite{Rol}. If $l=k$, $M(\sigma )$ is a
closed manifold; otherwise it has $k-l$ remaining cusps.

The Gromov--Thurston $2\pi $--theorem states that $M(\sigma )$ has a
complete metric of non--positive sectional curvature for 'most'
choices of $\sigma $. (Although the theorem was originally proved in
the context of $3$--manifolds, the same proof actually holds in any
dimension as observed in \cite{A}.) This means that the fundamental
group of $M(\sigma )$ is {\it semihyperbolic} in the sense of
\cite{AB}. The following immediate corollary of Theorem \ref{CEP}
shows that, in fact, $\pi _1(M(\sigma))$ is {\it hyperbolic relative
to finitely generated free abelian subgroups}, which is a much
stronger property than semihyperbolicity (see \cite{Reb}). We also
note that relatively hyperbolic groups of this type are Hopfian
\cite{Grov1,Grov2}, are $C^\ast $--exact \cite{Oz}, have finite
asymptotic dimension \cite{AsDim} (hence they satisfy the Novikov
Conjecture \cite{Yu}), and have many other nice properties. Below we
consider $\pi_1(T_i^{n-1})$ as a subgroup of $\pi _1(M)$, and set
$x_i=[\sigma _i]\in \pi_1(T_i^{n-1})$, $1\le i\le l$, and $x_i=1$
for $l<i\le k$.

\begin{cor}
There is a finite subset of nontrivial elements $\mathcal F\subset
\pi_{1}(M)$ such that if $x_i\notin \mathcal F$ for all $1\le i\le
k$, then the quotient groups $\pi_1(T_i^{n-1})/\langle x_i\rangle $
naturally inject into $M(\sigma )$ and $\pi _1(M(\sigma ))$ is
hyperbolic relative to the collection of finitely generated free
abelian subgroups $\{ \pi_1(T_i^{n-1})/\langle x_i\rangle \} _{1\le
i\le k}$.
\end{cor}

Let us discuss some algebraic applications of Theorem \ref{CEP}.
Recall that a subgroup $H$ of a group $G$ has the {\it Congruence
Extension Property} if for any $N\lhd H$, we have $H\cap N^G=N$ (or,
equivalently, the natural homomorphism $H/N\to G/N^G$ is injective).
An obvious example of the CEP is provided by the pair $G, H$, where
$H$ is a free factor of $G$. Another example is a cyclic subgroup
$H=\langle w\rangle $ generated by an arbitrary element $w$ of a
free group $F$. In this case the CEP for $H$ is equivalent to the
assertion that the element represented by $w$ has order $n$ in the
one relator group $F/\langle w^n\rangle ^F$, which is a part of the
well known theorem of Karrass, Magnus, and Solitar \cite{MKS}.
Olshanskii \cite{Ols95} noticed that the free group of rank $2$
contains subgroups of arbitrary rank having CEP. This easily implies
the Higman--Neumann-Neumann theorem stating that any countable group
can be embeded into a 2--generated group. The CEP has also been
extensively studied for semigroups and universal algebras (see
\cite{Blok,Tang,Wu} and references therein). It plays an important
role in some constructions of groups with 'exotic' properties
\cite{OS}.

We say that a subgroup $H$ of a group $G$ {\it almost has $CEP$} if
there is a finite set of non--trivial elements $\mathcal F\subseteq
H$ such that $H\cap N^G=N$ whenever $N\cap\mathcal F=\emptyset $.
Recall that a subgroup $H$ of a group $G$ is said to be {\it almost
malnormal}, if $H^g\cap H$ is finite for all $g\notin H$. Bowditch
\cite{Bow} proved that if $G$ is a hyperbolic group and $H$ is an
almost malnormal quasi--convex subgroup of $G$, then $G$ is
hyperbolic relative to $H$ (see also \cite{ESBG}). Thus the
following is an immediate corollary of Theorem \ref{CEP}.

\begin{cor}
Any almost malnormal quasi--convex subgroup of a hyperbolic group
almost has CEP.
\end{cor}

If $G$ is a free group, any almost malnormal subgroup $H\le G$ is
malnormal (i.e., it satisfies $H^g\cap H=\{ 1\} $ for all $g\notin
H$). It is also well--known that a subgroup of a finitely generated
free group is quasi--convex if and only if it is finitely generated.
Even the following result seems to be new.

\begin{cor}\label{CEPF}
Any finitely generated malnormal subgroup of a free group almost
has CEP.
\end{cor}

If the free group is finitely generated this is a particular case
of the previous corollary. To prove Corollary \ref{CEPF} in the
full generality, it suffices to notice that any finitely generated
subgroup $H$ of a free group $F$ belongs to a finitely generated
free factor $F_0$ of $F$ and $F_0$ has CEP as a subgroup of $F$.
This easily implies that  $H$ almost has CEP in $F$.

Considering a series of subgroups $K\lhd H\lhd F$ in a free group
$F$, where $K$ is not normal in $F$, it is easy to notice that the
word 'malnormal' can not be removed from the corollary. It is less
trivial that, in general, malnormal subgroups of free groups do not
have CEP. Here we sketch an example suggested by A. Klyachko. Let
$F$ be the free group with basis $x,y$. Using small cancellation
arguments, it is not hard to construct a malnormal subgroup $H$ of
$F$ generated by $x$ and some word $w\in [F,F]$. Then $\langle
x\rangle^F=H\ne\langle x\rangle ^H$ since $w\in \langle x\rangle
^F$.

Theorem \ref{CEP} also implies that, in an algebraic sense, the
group $G$ is approximated by its images obtained by peripheral
fillings. To be more precise, we recall that a group $G$ is {\it
fully residually} $\mathcal C$, where $\mathcal C$ is a class of
groups, if for any finite subset $S\subseteq G$, there is a
homomorphism of $G$ onto a group from $\mathcal C$ that is injective
on $S$. The study of this notion has a long history and is motivated
by the following observation: If $\mathcal C$ is a class of 'nice'
groups in a certain sense, then any (fully) residually $\mathcal C$
group also enjoys some nice properties.

Using Theorem \ref{CEP}, we will obtain some non--trivial examples
of fully residually hyperbolic groups. We recall that a group is
called {\it non--elementary} if it does not contain a cyclic
subgroup of finite index.

\begin{cor}\label{frh}
Suppose that a finitely generated group $G$ is hyperbolic relative
to a collection of subgroups $\Hl $ and for each $\lambda\in \Lambda
$, the group $H_\lambda $ is fully residually hyperbolic.  Then $G$
is fully residually hyperbolic. Moreover, if $G$ is non--elementary
and all subgroups $H_\lambda $ are proper, then $G$ is fully
residually non--elementary hyperbolic.
\end{cor}

For instance, fundamental groups of complete finite volume
Riemannian manifolds of pinched negative curvature are hyperbolic
relative to the cusp subgroups \cite{Bow,F}, which are virtually
nilpotent \cite{E}. It is well--know that any nilpotent group is
residually finite \cite{H} and hence so is any virtually nilpotent
group. Finally we recall that finite groups are hyperbolic.
Combining this with Corollary \ref{frh} we obtain
\begin{cor}
Fundamental groups of complete finite volume Riemannian manifolds of
pinched negative curvature are fully residually non--elementary
hyperbolic.
\end{cor}

Note that any fully residually non--elementary hyperbolic group $G$
has infinite quotients of bounded period and, moreover,
$\bigcap\limits _{n=1}^\infty G^n=\{1 \} $, where $G^n=\{ g^n\, | \,
g\in G\} $. This easily follows from the result of Ivanov and
Olhanskii \cite{IO}.

Another application is related to the well known question of whether
all hyperbolic groups are residually finite. Although in many
particular cases the answer is known to be positive (see \cite{W}
and references therein), in the general case the question is still
open. The following obvious consequence of Corollary \ref{frh} shows
that this problem is equivalent to its relative analogue. In
particular, in order to construct a non--residually finite
hyperbolic group it suffices to find a non--residually finite group
that is hyperbolic relative to a collection of residually finite
subgroups.

\begin{cor}
The following assertions are equivalent.
\begin{enumerate}

\item Suppose that a finitely generated group $G$ is hyperbolic
relative to a collection of residually finite subgroups. Then $G$
is residually finite.

\item Any hyperbolic group is residually finite.
\end{enumerate}
\end{cor}

The paper is organized as follows. In the next section we give the
definition of relatively hyperbolic groups and provide a background
for the rest of the paper. The proof of the main theorem consists of
two ingredients. The first one is Proposition \ref{s(n)} concerning
geodesic polygons in Cayley graphs of relatively hyperbolic groups.
It is proved in Section 3 and seems to be of independent interest.
The second ingredient is a surgery on van Kampen diagrams described
in Sections 4. Theorem \ref{CEP} and Corollary \ref{frh} are proved
in Section 5.

\bigskip

{\bf Acknowledgments.} The author is grateful to Koji Fujiwara,
Anton Klyachko, Alexander Olshanskii, and John Ratcliffe for useful
discussions.

%%%%%%%%%%%%%%%%%%%%%%%%%%%%%%%%%%%%%%%%%%%%%%%%%%%%%%%%%%%%%%%

\section{Preliminaries}

%%%%%%%%%%%%%%%%%%%%%%%%%%%%%%%%%%%%%%%%%%%%%%%%%%%%%%%%%%%%%%%

{\bf Some conventions and notation.} Given a word $W$ in an
alphabet $\mathcal A$, we denote by $\| W\| $ its length. We also
write $W\equiv V$ to express the letter--for--letter equality of
words $W$ and $V$. Recall that a subset $X$ of a group $G$ is said
to be {\it symmetric} if for any $x\in X$, we have $x^{-1}\in X$.
In this paper all generating sets of groups under consideration
are supposed to be symmetric.

\vspace{3mm}

\noindent {\bf Word metrics and Cayley graphs.} Let $G$ be a group
generated by a (symmetric) set $\mathcal A$. Recall that the {\it
Cayley graph} $\Ga $ of a group $G$ with respect to the set of
generators $\mathcal A$ is an oriented labelled 1--complex with
the vertex set $V(\Ga )=G$ and the edge set $E(\Ga )=G\times
\mathcal A$. An edge $e=(g,a)$ goes from the vertex $g$ to the
vertex $ga$ and has label $\phi (e)\equiv a$. As usual, we denote
the origin and the terminus of the edge $e$ by $e_-$ and $e_+$
respectively. Given a combinatorial path $p=e_1e_2\ldots e_k$ in
the Cayley graph $\Ga $, where $e_1, e_2, \ldots , e_k\in E(\Ga
)$, we denote by $\phi (p)$ its label. By definition, $\phi
(p)\equiv \phi (e_1)\phi(e_2)\ldots \phi (e_k).$ We also denote by
$p_-=(e_1)_-$ and $p_+=(e_k)_+$ the origin and the terminus of $p$
respectively. The length $l(p)$ of $p$ is the number of edges in
$p$.

The {\it (word) length} $|g|_\mathcal A$ of an element $g\in G$ is
defined to be the length of a shortest word in $\mathcal A$
representing $g$ in $G$. This defines a metric on $G$ by
$dist_\mathcal A(f,g)=|f^{-1}g|_\mathcal A$. We also denote by
$dist _\mathcal A$ the natural extension of the corresponding
metric on the Cayley graph $\Ga $.

\vspace{3mm}

\noindent {\bf Van Kampen Diagrams.} For technical reasons, it is
convenient to define the diagrams as in the book \cite{Ols-book},
i.e., to allow the so called $0$--cells. More precisely, a {\it
van Kampen diagram} $\Delta $ over a presentation
\begin{equation}
G=\langle \mathcal A\; | \; \mathcal O\rangle \label{ZP}
\end{equation}
is a finite oriented connected planar 2--complex endowed with a
labelling function $\phi : E(\Delta )\to \mathcal A\sqcup \{ 1\} $,
where $E(\Delta ) $ denotes the set of oriented edges of $\Delta $,
such that $\phi (e^{-1})\equiv (\phi (e))^{-1}$. The symbol $1$
denotes the trivial word here. We call the edges of $\Delta $
labelled by letters from $\mathcal A$ {\it essential}. Labels of
paths are defined as in the case of Cayley graphs, but the symbols
$1$ are always omitted. Hence labels of paths are words in $\mathcal
A$. When defining the lengths of a paths, we do not count the edges
labelled by $1$. Thus the lengths of a path always agrees with the
lengths of its label.

Given a cell $\Pi $ of $\Delta $, we denote by $\partial \Pi$ the
boundary of $\Pi $; similarly, $\partial \Delta $ denotes the
boundary of $\Delta $. The labels of $\partial \Pi $ and $\partial
\Delta $ are defined up to cyclic permutations. An additional
requirement is that any cell $\Pi $ of $\Delta $ satisfies one of
the following two conditions:
\begin{enumerate}
\item $\phi (\partial \Pi )$ is equal to (a cyclic permutation of)
a word $P^{\pm 1}$, where $P\in \mathcal O$. We call such cells
{\it essential}.

\item  The set of essential edges of $\partial \Pi $ is empty or
consists of exactly two edges whose labels are $a$ and $a^{-1}$ for
some $a\in \mathcal A$. Thus $\phi (\Pi )$ represents the identity
element in the free group generated by $\mathcal A$. The cells of
this type are called {\it $0$--cells}.
\end{enumerate}

\begin{figure}
 \vspace{2mm}\hspace{15mm} \includegraphics{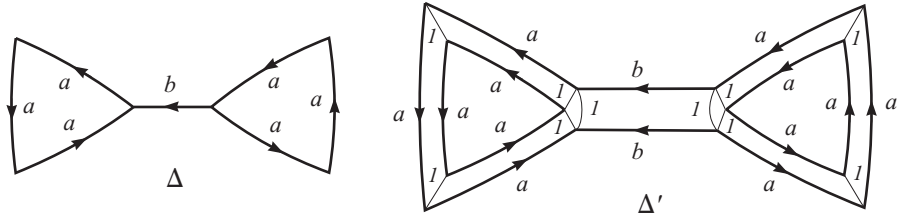}\\
\vspace{-2mm}
 \caption{A $0$--refinement of a diagram
over the presentation $\langle a,b\; |\; a^3=1\rangle $.}
\label{figa}
\end{figure}

One says that a diagram $\Delta ^\prime $ is a {\it
$0$--refinement} of a diagram $\Delta $, if, roughly speaking, it
is obtained from $\Delta $ by replacing some edges and vertices of
$\Delta $ with appropriate $0$--cells (see Fig. \ref{figa}).  This
notion is quite standard and we refer the reader to \cite[Ch.
4]{Ols-book} for details.

The van Kampen Lemma states that a word $W$ over an alphabet
$\mathcal A$ represents the identity in the group given by
(\ref{ZP}) if and only if there exists a connected
simply--connected planar diagram $\Delta $ over (\ref{ZP}) such
that $\phi (\partial \Delta )\equiv W$ \cite[Ch. 5, Theorem
1.1]{LS}. Using $0$--refinement if necessary, we may always assume
$\Delta $ to be homeomorphic to a disk.

Finally we make the following quite obvious observation (see
\cite[Ch. 4]{Ols-book}).

\begin{lem}\label{map}
Let $\Delta $ be a connected simply--connected diagram over
(\ref{ZP}) with a fixed vertex $O$, $\Ga $ the Cayley graph of $G$
with respect to the generating set $\mathcal A$. Then there is a
(unique) continuous map $\mu \colon Sk^{(1)}\, (\Delta )\to \Ga $
that maps $O$ to the identity vertex of $\Ga $, vertices of $\Delta
$ and edges labelled by $1$ to vertices of $\Ga $, and essential
edges of $\Delta $ to edges of $\Ga $ preserving labels and
orientation.
\end{lem}

\vspace{3mm}

\noindent {\bf Hyperbolic spaces.} Recall that a metric space $X$
is {\it $\delta $--hyperbolic} for some $\delta \ge 0$ (or simply
{\it hyperbolic}) if for any geodesic triangle $T$ in $X$, any
side of $T$ belongs to the union of the closed $\delta
$--neighborhoods of the other two sides \cite{Gro}.

In this paper we use some results about polygons in hyperbolic
spaces. We recall that a path $p$ in a metric space $X$ is called
{\it $(\lambda , c)$--quasi--geodesic} for some $\lambda \ge 1$,
$c\ge 0$ if
$$l(q)\le  \lambda dist(q_-, q_+)+c$$ for any subpath $q$ of $p$.
Recall also that for any fixed $\delta $, $\lambda $, and $c$, all
$(\lambda , c)$--quasi--geodesics with same endpoints in a $\delta
$--hyperbolic space are uniformly close (see, for example,
\cite[Ch. III. H, Theorem 1.7]{BH}). The first lemma can easily be
derived from this result and the definition of a hyperbolic space
by drawing the diagonal.

\begin{lem}\label{rect}
For any $\delta \ge 0$, $\lambda \ge 1$, $c\ge 0$, there exists a
constant $\kappa =\kappa(\delta , \lambda , c)\ge 0$ with the
following property. Let $Q$ be a quadrangle in a
$\delta$--hyperbolic space whose sides are $(\lambda ,
c)$--quasi--geodesic. Then each side of $Q$ belongs to the closed
$\kappa $--neighborhood of the union of the other three sides.
\end{lem}

The next lemma was proved by Olshanskii \cite[Lemma 23]{Ols92} for
geodesic polygons. In \cite{Ols92}, the inequality (\ref{distuv})
had the form $dist (u,v)\le 2\delta ( 2+\log_2 n)$. Passing to
quasi--geodesic polygons we only need to add a constant to the right
hand side according to the above--mentioned property of
quasi--geodesics in hyperbolic spaces.

\begin{lem}\label{P}
For any $\delta \ge 0$, $\lambda \ge 1$, $c\ge 0$, there exists a
constant $\theta =\theta (\delta , \lambda , c)$ with the following
property. Let $\P =p_1\ldots p_n$ be a $(\lambda ,
c)$--quasi--geodesic $n$--gon in a $\delta $--hyperbolic space. Then
there are points $u$ and $v$ on sides of $\mathcal P$ such that
\begin{equation}\label{distuv}
dist (u,v)\le 2\delta ( 2+\log_2 n)+\theta
\end{equation}
and the geodesic segment connecting $u$ to $v$ divides $\P $ into an
$m_1$--gon and $m_2$--gon such that $n/4 < m_i < 3n/4 +2$.
\end{lem}

\vspace{3mm}

\noindent {\bf Relatively hyperbolic groups.} In this paper we use
the notion of relative hyperbolicity whic is sometimes called strong
relative hyperbolicity and goes back to Gromov \cite{Gro}. There are
many equivalent definitions of (strongly) relatively hyperbolic
groups \cite{Bow,DS,F,RHG}. We recall the isoperimetric
characterization suggested in \cite{RHG}, which is most suitable for
our purposes.

Let $G$ be a group, $\Hl $ a collection of subgroups of $G$, $X$ a
subset of $G$. We say that $X$ is a {\it relative generating set of
$G$ with respect to $\Hl $} if $G$ is generated by $X$ together with
the union of all $H_\lambda $. (In what follows we always assume $X$
to be symmetric.) In this situation the group $G$ can be regarded as
a quotient group of the free product
\begin{equation}
F=\left( \ast _{\lambda\in \Lambda } H_\lambda  \right) \ast F(X),
\label{F}
\end{equation}
where $F(X)$ is the free group with the basis $X$. If the kernel of
the natural homomorphism $F\to G$ is a normal closure of a subset
$\mathcal R$ in the group $F$, we say that $G$ has {\it relative
presentation}
\begin{equation}\label{G}
\langle X,\; H_\lambda, \lambda\in \Lambda \; |\; \mathcal R
\rangle .
\end{equation}
If $\sharp\, X<\infty $ and $\sharp\, \mathcal R<\infty $, the
relative presentation (\ref{G}) is said to be {\it finite} and the
group $G$ is said to be {\it finitely presented relative to the
collection of subgroups $\Hl $.}

Set
\begin{equation}\label{H}
\mathcal H=\bigsqcup\limits_{\lambda\in \Lambda} (H_\lambda
\setminus \{ 1\} ) .
\end{equation}
Given a word $W$ in the alphabet $X\cup \mathcal H$ such that $W$
represents $1$ in $G$, there exists an expression
\begin{equation}
W=_F\prod\limits_{i=1}^k f_i^{-1}R_i^{\pm 1}f_i \label{prod}
\end{equation}
with the equality in the group $F$, where $R_i\in \mathcal R$ and
$f_i\in F $ for $i=1, \ldots , k$. The smallest possible number
$k$ in a representation of the form (\ref{prod}) is called the
{\it relative area} of $W$ and is denoted by $Area^{rel}(W)$.

\begin{defn}
A group $G$ is {\it hyperbolic relative to a collection of
subgroups} $\Hl $ if $G$ is finitely presented relative to $\Hl $
and there is a constant $C>0$ such that for any word $W$ in $X\cup
\mathcal H$ representing the identity in $G$, we have $Area^{rel}
(W)\le C\| W\| $. The constant $C$ is called an {\it isoperimetric
constant} of the relative presentation (\ref{G}).
\end{defn}

In particular, $G$ is an ordinary {\it hyperbolic group} if $G$ is
hyperbolic relative to the trivial subgroup. An equivalent
definition says that $G$ is hyperbolic if it is generated by a
finite set $X$ and the Cayley graph $\Gamma (G, X)$ is hyperbolic.
In the relative case these approaches are not equivalent, but we
still have the following \cite[Theorem 1.7]{RHG}.

\begin{lem}\label{CG}
Suppose that $G$ is a group hyperbolic relative to a collection of
subgroups $\Hl $. Let $X$ be a finite relative generating set of
$G$ with respect to $\Hl $. Then the Cayley graph $\G $ of $G$
with respect to the generating set $X\cup \mathcal H$ is a
hyperbolic metric space.
\end{lem}

Observe also that the relative area of a word $W$ representing $1$
in $G$ can be defined geometrically via van Kampen diagrams. Let
$G$ be a group given by the relative presentation (\ref{G}) with
respect to a collection of subgroups $\Hl $.  We denote by
$\mathcal S$ the set of all words in the alphabet $\mathcal H$
representing the identity in the groups $F$ defined by (\ref{F}).
Then $G$ has the ordinary (non--relative) presentation
\begin{equation}\label{Gfull}
G=\langle X\cup\mathcal H\; |\;\mathcal S\cup \mathcal R \rangle .
\end{equation}
A cell in van Kampen diagram $\Delta $ over (\ref{Gfull}) is called
an {\it $\mathcal R$--cell} if its boundary is labeled by a word
from $\mathcal R$. We denote by $N_\mathcal R(\Delta )$ the number
of $\mathcal R$--cells of $\Delta $. Obviously given a word $W$ in
$X\cup\mathcal H$ that represents $1$ in $G$, we have
$$
Area^{rel}(W)=\min\limits_{\phi (\partial \Delta ) \equiv W}
N_\mathcal R (\Delta ) ,
$$
where the minimum is taken over all disk van Kampen diagrams with
boundary label $W$.

We recall an auxiliary terminology introduced in \cite{RHG}, which
plays an important role in our paper.

\begin{defn}
Let $q$ be a path in the Cayley graph $\G $. A (non--trivial)
subpath $p$ of $q$ is called an {\it $H_\lambda $--subpath} for some
$\lambda \in \Lambda $, if the label of $p$ is a word in the
alphabet $H_\lambda\setminus \{ 1\} $. If $p$ is a maximal
$H_\lambda $--subpath of $q$, i.e. it is not contained in a bigger
$H_\lambda $--subpath, then $p$ is called an {\it $H_\lambda
$--component} (or simply a {\it component}) of $q$.

Two $H_\lambda $--subpaths (or $H_\lambda $--components) $p_1, p_2$
of a path $q$ in $\G $ are called {\it connected} if there exists a
path $c$ in $\G $ that connects some vertex of $p_1$ to some vertex
of $p_2$ and ${\phi (c)}$ is a word consisting of letters from $
H_\lambda\setminus\{ 1\} $. In algebraic terms this means that all
vertices of $p_1$ and $p_2$ belong to the same coset $gH_\lambda $
for a certain $g\in G$. Note that we can always assume that $c$ has
length at most $1$, as every nontrivial element of $H_\lambda
\setminus\{ 1\} $ is included in the set of generators.  An
$H_\lambda $--component $p$ of a path $q$ is called {\it isolated }
(in $q$) if no distinct $H_\lambda $--component of $q$ is connected
to $p$.
\end{defn}

To every subset $\Omega$ of $G$, we can associate a (partial)
distance function $dist_\Omega \colon G\times G\to [0, \infty ]$
as follows. If $g_1, g_2\in G$ and $g_1^{-1}g_2\in \langle \Omega
\rangle $, we define $dist _\Omega (g_1, g_2)=|g_1^{-1}g_2|_\Omega
$, where $|\cdot |_\Omega $ is the word length with respect to
$\Omega $. If $g_1^{-1}g_2\notin \langle \Omega \rangle $, we set
$dist _\Omega (g_1, g_2) =\infty $. Finally, for any path $p$ in
$\G $, we define its $\Omega $--length as $$l_\Omega (p)=dist
_\Omega (p_-, p_+).$$

The lemma below was proved in \cite[Lemma 2.27]{RHG} in the case
when $p_1, \ldots , p_k$ are $H_\lambda $--components for a fixed
$\lambda $. Actually the proof from \cite{RHG} works in the
general case as well. Here we provide it for convenience of the
reader.

\begin{lem}\label{Omega}
Let $G$ be a group that is hyperbolic relative to a collection of
subgroups $\Hl $. Then there exists a finite subset $\Omega
\subseteq G$ and a constant $L>0$ such that the following condition
holds. Let $q$ be a cycle in $\G $, $p_1, \ldots , p_k$ a set of
isolated components of $q$. Then the $\Omega $--lengths of $p_i$'s
satisfy
$$ \sum\limits_{i=1}^k l_\Omega (p_i)\le Ll(q).$$
\end{lem}

\begin{figure}
 \hspace{36mm} \includegraphics{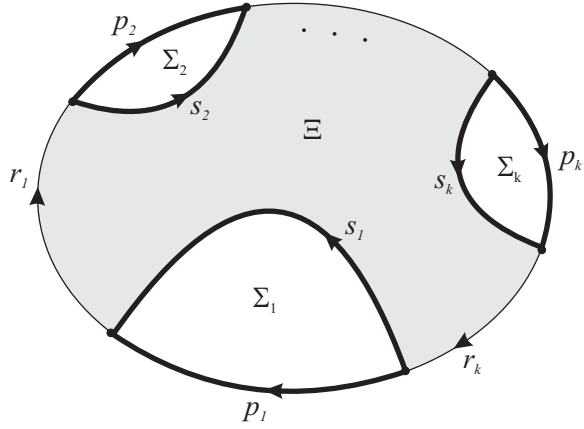}\\
\vspace{-5mm} \caption{Decomposition of the diagram $\Delta
$}\label{fig0}
\end{figure}

\begin{proof}
Let $\Omega $ consist of all letters from $\mathcal H$ that appear
in words from the set $\mathcal R$ (see (\ref{G})). Since $\sharp
\, \mathcal R<\infty $, we have $\sharp\, \Omega <\infty $. To
prove the lemma we consider a van Kampen diagram $\Delta $ over
(\ref{Gfull}) whose boundary label is $\phi (q)$. In what follows
we identify $\partial \Delta $ with $q$.

Assume that $q=p_1r_1\cdots p_kr_k$ and $p_i$ is an $H_{\lambda
_i}$--component for some $\lambda _i\in \Lambda $. For each $i=1,
\ldots , k$, let $\mathcal D_i$ denote the set of all subdiagrams
of $\Delta $ bounded by $p_i(p^{\prime }_i)^{-1}$, where
$p^{\prime }_i$ is a simple path in $\Delta $ such that
$(p^{\prime }_i)_-=(p_i)_-$, $(p^{\prime }_i)_+=(p_i)_+$, and
$\phi (p^{\prime }_i)$ is a word in the alphabet $H_{\lambda
_i}\setminus \{ 1\} $. We choose a subdiagram $\Sigma _i \in
\mathcal D_i$ that has maximal number of cells among all
subdiagrams from $\mathcal D_i$.

Let $\partial \Sigma _i= p_is_i^{-1}$. Since $p_i$ is an isolated
component of $q$, the path $s_i$ has no common edges with $r_i$,
$i=1, \ldots k$, and the sets of edges of $s_i$ and $s_j$ are
disjoint whenever $j\ne i$. Therefore each edge $e$ of $s_i$ belongs
to a boundary of some cell $\Pi $ of the subdiagram $\Xi $ of
$\Delta $ bounded by $s_1r_1\cdots s_kr_k $. If $\Pi $ is an
$S$--cell, then $\phi (\Pi )$ is a word in the alphabet $H_{\lambda
_i}\setminus \{ 1\} $. Hence by joining $\Pi $ to $\Sigma _i$ we get
a subdiagram $\Sigma _i^\prime \in \mathcal D_i$ with bigger number
of cells that contradicts the choice of $\Sigma _i $. Thus each edge
of $s_i$ belongs to a boundary of an $\mathcal R$--cell and, in
particular, has $\Omega $--length $1$. The  total number of such
edges does not exceed the number of $\mathcal R$--cells in $\Xi $
times the maximal number of edges in boundary of an $\mathcal
R$--cell. Therefore we have
$$ \sum\limits_{i=1}^k l_\Omega (p_i)= \sum\limits_{i=1}^k
l_\Omega (s_i)\le M Area ^{rel} (\phi(\partial \Delta ))\le
MCl(q),$$ where $C$ is the isoperimetric constant of (\ref{G}) and
$M=\max\limits_{R\in \mathcal R} \| R\| .$
\end{proof}

%%%%%%%%%%%%%%%%%%%%%%%%%%%%%%%%%%%%%%%%%%%%%%%%%%%%%%%%%%%%%%%%

\section{Components and quasi--geodesic polygons.}

%%%%%%%%%%%%%%%%%%%%%%%%%%%%%%%%%%%%%%%%%%%%%%%%%%%%%%%%%%%%%%%%

Throughout the rest of the paper let $G$ denote a group that is
hyperbolic relative to a collection of subgroups $\Hl $. Let also
$X$ be a finite generating set of $G$ with respect to $\Hl $,
$\Omega $ the subset of $G$, and $L$ the constant provided by Lemma
\ref{Omega}. In this section we show that the bound on $\Omega
$--lengths of components in Lemma \ref{Omega} can be essentially
improved in some special cases.

\begin{defn}
For $\lambda \ge 1$, $c\ge 0$, and $n\ge 2$, let $\mathcal
Q_{\lambda , c}(n)$ denote the set of all pairs $(\P ,\, I)$,
where $\P =p_1\ldots p_n$ is an $n$--gon in $\G $ and $I$ is a
distinguished subset of the set of sides $\{ p_1, \ldots, p_n\} $
of $\P $ such that:
\begin{enumerate}
\item Each side $p_i\in I$ is an isolated component of $\P $.

\item Each side $p_i\notin I$ is $(\lambda ,c)$--quasi--geodesic.
\end{enumerate}
For technical reasons, it is convenient to allow some of the sides
$p_1, \ldots , p_n$ to be trivial. Thus we have $\mathcal
Q_{\lambda , c}(2)\subseteq\mathcal Q_{\lambda , c}(3)\subseteq
\ldots $. Below we also use the following notation for vertices of
$\P $:
$$ x_1=(p_n)_+=(p_1)_-, \; x_2=(p_1)_+=(p_2)_-,\; \ldots , \;
x_n=(p_{n-1})_+=(p_{n})_-.$$ Given $(\P , I)\in \mathcal
Q_{\lambda , c}(n)$, we set
$$s(\P , I)=\sum\limits_{p_i\in I} l_\Omega (p_i)$$ and consider
the quantity $$ s_{\lambda , c}(n)=\sup\limits_{(\P , I) \in
\mathcal Q_{\lambda , c}(n)} s(\P , I).$$
\end{defn}

Observe that, a priori, it is not clear whether $s_{\lambda , c}(n)$
is finite for fixed values of $n$, $\lambda $, and $c$. The main
purpose of this section is to prove a much stronger result.

\begin{prop}\label{s(n)}
For any $\lambda \ge 1$, $c\ge 0$, there exists a constant
$D=D(\lambda , c)>0$ such that $s_{\lambda , c}(n)\le Dn $ for any
$n\in \mathbb N$.
\end{prop}

The following simple observation will often be used in this
section without special references. If $p_1$, $p_2$ are connected
components of some path in $\G $, then for any two vertices $u\in
p_1$ and $v\in p_2$, we have $\dxh (u,v)\le 1$ .

The proof of Proposition \ref{s(n)} is by induction on $n$. We
begin with the case $n\le 4$.

\begin{lem}\label{l1}
For any $\lambda \ge 1$, $c\ge 0$, and $n\le 4$,  $s_{\lambda ,
c}(n)$ is finite.
\end{lem}

\begin{proof}
Suppose that $(\P ,\, I)\in \mathcal Q_{\lambda , c}(4)$, $\P
=p_1p_2p_3p_4$. We want to show that $s(\P , I)$ is bounded by a
constant, which depends on $\lambda $, $c$, and the hyperbolicity
constant $\delta $ of the graph $\G $ only. We notice that for a
path $p$ in $\G$, $l_\Omega (p)$ depends on $p_-$ and $p_+$ only
and otherwise is independent of $p$ itself. Since every $p_i\in I$
is an $H_{\lambda _i}$--component for some $\lambda _i\in \Lambda
$, we can replace each $p_i\in I$ with a single edge $e_i$
labelled by an appropriate element of $H_{\lambda _i}$. Clearly
$e_i$ is isolated in $\P $ whenever $p_i$ is. Thus we may assume
that $l(p_i)=1$ whenever $p_i\in I$.

Let $\kappa =\kappa (\delta , \lambda , c) $ be the constant
provided by Lemma \ref{rect}. Without loss of generality we may
assume $\kappa $ to be a positive integer. According to Lemma
\ref{Omega}, it suffices to show that for each $p_i\in I$, there
is a cycle $c_i$ in $\G $ of length less than $K=100(\lambda
\kappa +c+\kappa )$ such that $p_i$ is an isolated component of
$c_i$. There are $4$ cases to consider.

\medskip

\noindent {\bf Case 1.}  Suppose $\sharp\, I=4$. Then the
assertion of the lemma is obvious. Indeed $l(\P )=4<K$ as each
$p_i\in I$ has lengths $1$, and we can set $c_i=\P $ for all $i$.

\medskip

\noindent {\bf Case 2.}  Suppose $\sharp\, I=3$, say $I=\{
p_1,p_2,p_3\} $. Since $p_4$ is $(\lambda, c)$--quasi--geodesic,
we have
$$l(p_4)\le \lambda \dxh (x_4, x_1)+c\le 3\lambda +c$$ by
the triangle inequality. Hence $l(\P )\le 3\lambda +c+3< K$ and we
can set $c_i=\P $ again.

\medskip

\noindent {\bf Case 3.}  Assume now that $\sharp\, I=2$. Up to
enumeration of the sides, there are two possibilities to consider.

\smallskip

a) First suppose $I=\{ p_1,p_2\} $. If $\dxh (x_3, x_4)< \kappa
+2$, we have
$$
l(p_3)\le \lambda \dxh (x_3, x_4) +c<\lambda (\kappa +2)+c,
$$
$$
\begin{array}{cl}
l(p_4)\le & \lambda \dxh (x_4, x_1) +c \le \\ &  \lambda (\dxh
(x_1, x_2) +\dxh (x_2, x_3) +\dxh (x_3, x_4))+ c< \\ & \lambda (1
+1 + \kappa +2) +c \le \lambda (\kappa +4)+c,
\end{array}
$$
and hence $$l(\P )<2+ l(p_3)+l(p_3) <\lambda (2\kappa +6)+2c+2<
K.$$ Thus we may assume $\dxh (x_3, x_4)\ge \kappa +2$. Let $u$ be
a vertex on $p_3$ such that $\dxh (x_3, u)=\kappa +2$. By Lemma
\ref{rect} there exists a vertex $v\in p_1\cup p_2\cup p_4 $ such
that $\dxh (u,v)\le \kappa $. Note that, if fact, $v\in p_4$.
Indeed otherwise $v=x_2$ or $v=x_3$ and we have
$$\dxh (x_3, u)\le \dxh (x_3, v)+\dxh (u, v)\le 1+\kappa $$ that
contradicts the choice of $u$.

\begin{figure}
\hspace{31mm} \includegraphics{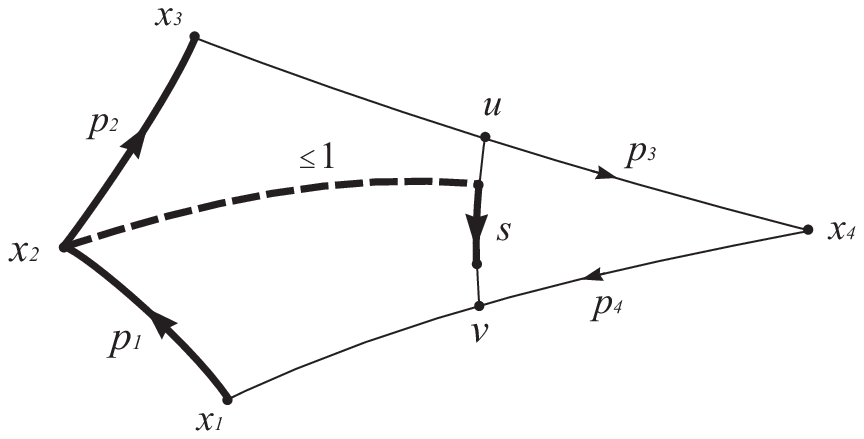} \\
\hspace*{36mm} \includegraphics{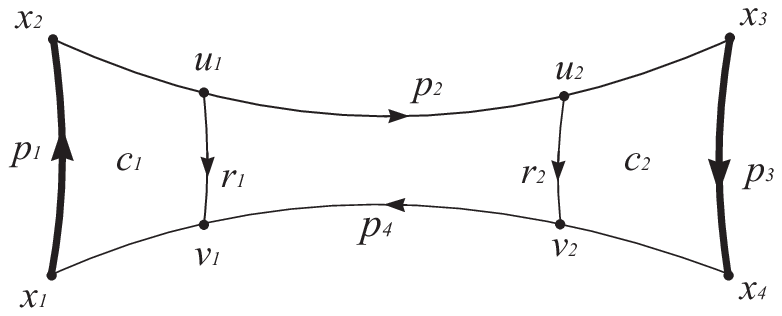}\\
\vspace*{-8mm} \caption{Cases 3 a) and b)}\label{fig1}
\end{figure}

Let $r$ be a geodesic path in $\G $ connecting $u$ to $v$. We wish
to show that no component of $r$ is connected to $p_1$ or $p_2$.
Indeed suppose that a component $s$ of $r$ is connected to $p_1$
or $p_2$ (Fig.\ref{fig1}). Then $\dxh (x_2, s_-)\le 1$ and we
obtain
$$
\begin{array}{rl}
\dxh (u, x_3)\le & \dxh (u, s_-)+\dxh (s_-, x_2)+ \dxh (x_2,
x_3)\le \\ & (\kappa -1)+1 +1=\kappa +1.
\end{array}
$$
This contradicts the choice of $u$ again. Note also that $p_1$,
$p_2$ can not be connected to a component of $p_3$ or $p_4$ as
$p_1$, $p_2$ are isolated components in $\P $. Therefore $p_1$ and
$p_2$ are isolated components of the cycle
$$c=p_1p_2[x_3, u]r[v, x_1],$$ where $[x_3, u]$ and $[v, x_1]$ are
segments of $p_3$ and $p_4$ respectively. Using the triangle
inequality, it is easy to check that $l([v, x_1]) \le
\lambda(2\kappa +4)$ and $l(c)\le \lambda (3\kappa +6)+2c+\kappa
+2< K$.

\smallskip

 b) Let $I=\{ p_1,p_3\}  $. If $\dxh (x_2, x_3)< 2\kappa +2$,
we obtain $l(\P )< K$ arguing as in the previous case. Now assume
that $\dxh (x_2, x_3)\ge 2\kappa +2$. Let $u_1$ (respectively
$u_2$) be the vertex on $p_2$ such that $\dxh (x_2, u_1)=\kappa
+1$ (respectively $\dxh (x_3, u_2)=\kappa +1$). By Lemma
\ref{rect} there exist vertices $v_1, v_2$ on $p_1\cup p_3\cup
p_4$ such that $\dxh (v_i,u_i)\le \kappa $, $i=1,2$. In fact,
$v_1, v_2$ belong to $p_4$ (Fig.\ref{fig1}). Indeed the reader can
easily check that the assumption $v_1=x_2$ (respectively
$v_1=x_3$) leads to the inequality $\dxh (x_2, u_1)\le \kappa $
(respectively $\dxh (x_2, x_3)\le 2\kappa +1$). In both cases we
get a contradiction. Hence $v_1\in p_4$ and similarly $v_2\in
p_4$.

Let $r_i$, $i=1,2$, be a geodesic path in $\G $ connecting $u_i$
to $v_i$. We set
$$c_1=p_1[x_2,u_1]r_1[v_1, x_1]$$ and
$$c_3=p_3[x_4, v_2]r_2^{-1}[u_2, x_3].$$
Arguing as in Case 3a) we can easily show that $p_i$ is an
isolated component of $c_i$ and $l(c_i)<K$ for $i=1,2$.

\medskip

\begin{figure}
   \hspace{12mm} \includegraphics{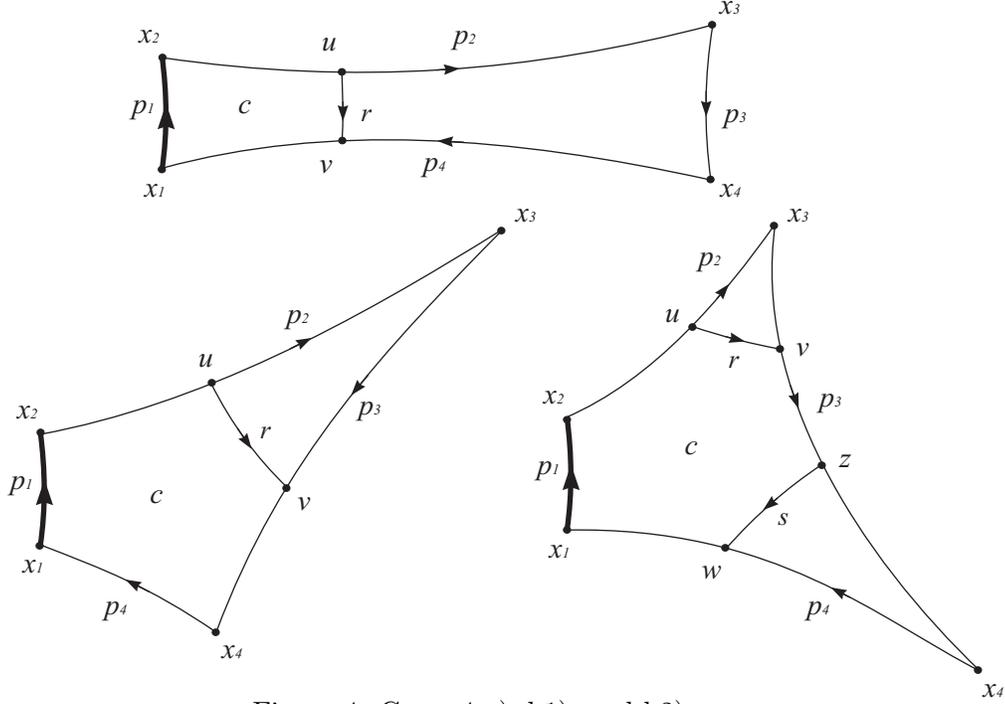}\\
  \vspace*{-10mm}
  \caption{Cases 4 a),  b1), and b2). }\label{fig3}
\end{figure}

\noindent {\bf Case 4.}  Finally assume $\sharp I=1$. To be
definite, let $I=\{ p_1\} $. If $\dxh(x_2,x_3)< \kappa +1$ and
$\dxh(x_4,x_1)< \kappa +1$, we obtain $l(\P )< K$ as in the
previous cases. Thus, changing the enumeration of the sides if
necessary, we may assume that $\dxh(x_2,x_3)\ge \kappa +1$. Let
$u$ be a point on $p_2$ such that $\dxh (x_2, u)=\kappa +1$, $v$ a
point on $p_1\cup p_3\cup p_4$ such that $\dxh (u,v)\le \kappa $,
$r$ a geodesic path in $\G $ connecting $u$ to $v$. As above it is
easy to show that $v\in p_3\cup p_4$. Let us consider two
possibilities (see Fig. \ref{fig3}).

\smallskip

a) $v\in p_4$. Using the same arguments as in Cases 2 and 3 the
reader can easily prove that $p_1$ is an isolated component of the
cycle
\begin{equation}\label{c4a}
c=p_1[x_2,u]r[v, x_1].
\end{equation}
It is easy to show that $l(c)<K$.

\smallskip

b) $v\in p_3$. Here there are 2 cases again.

\smallskip

b1) If $\dxh (x_1,x_4)< \kappa +1$, then we set
$$
c=p_1[x_2,u]r[v,x_4]p_4.
$$
The standard arguments show that $l(c)<K$ and $p_1$ is isolated in
$c$.

\smallskip

b2) $\dxh (x_1,x_4)\ge \kappa  +1$.  Let $w$ be a vertex on $p_4$
such that $\dxh (x_1, w)=\kappa +1$, $z$ a vertex on $p_1\cup
p_2\cup p_3$ such that $\dxh (z,w)\le \kappa $. Again, in fact,
our assumptions imply that $z\in p_2\cup p_3$. If $z\in p_2$, the
lemma can be proved by repeating the arguments from the case 4a)
(after changing enumeration of the sides). If $z\in p_3$, we set
$$c=p_1[x_2,u]r[v,z]s[w,x_1],$$ where $s$ is a geodesic in $\G $
connecting $z$ to $w$. It is straightforward to check that $p_1$
is an isolated component of $c$ and $l(c)<K$. We leave details to
the reader.

\end{proof}

\begin{lem}\label{l2}
For any $n\ge 4 $, we have
\begin{equation}\label{ge5}
s_{\lambda , c}(n)\le n (s_{\lambda , c}(n-1)+s_{\lambda , c}(4)).
\end{equation}
\end{lem}

\begin{figure}
  \hspace{37mm}\includegraphics{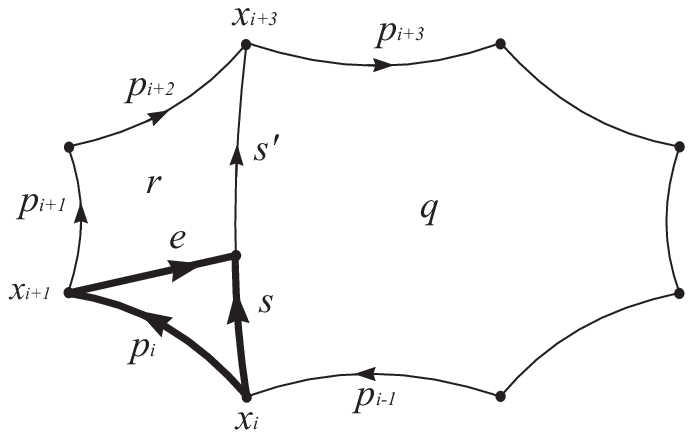}
  \caption{}\label{fig4}
\end{figure}

\begin{proof}
We proceed by induction on $n$. The case $n=4$ is obvious, so we
assume that $n\ge 5$. Let $(\P , I)\in \mathcal Q_{\lambda ,
c}(n)$, $p_i\in I$, and let $q$ be a geodesic in $\G $ connecting
$x_i$ to $x_{i+3}$ (indices are taken $mod\, n$). If $p_i$ is
isolated in the cycle $p_ip_{i+1}p_{i+2}q^{-1}$, we have $l_\Omega
(p_i)\le s_{\lambda , c}(4)$. Assume now that the component $p_i$
is not isolated in the cycle $p_ip_{i+1}p_{i+2}q^{-1}$. As $p_i$
is isolated in $\P $, this means that $p_i$ is connected to a
component $s$ of $q$. Hence $\dxh (x_i, s_+)\le 1$. Since $q$ is
geodesic in $\G $, this implies $s_-=x_i$ (see Fig. \ref{fig4}).

Let $q=ss^\prime $ and let $e$ denote a paths in $\G $ of lengths
at most $1$ such that $e_-=x_{i+1}$, $e_+=s_+$, and $\phi (e)$ is
a word in $\mathcal H$. We notice that if $e$ is nontrivial, it is
an isolated component of the cycle
$r=p_{i+1}p_{i+2}(s^{\prime})^{-1}$. Indeed if $e$ is connected to
a component of $p_{i+1}$ or $p_{i+2}$, then $p_i$ is not isolated
in $p$, and if $e$ is connected to a component of $s^\prime $,
then $q$ is not geodesic.  Similarly $s$ is an isolated component
of $p_{i+3}\ldots p_{i-1} ss^\prime $. Hence $l_\Omega (s)\le
s_{\lambda , c}(n-1)$ by the inductive assumption and $l_\Omega
(e)\le s_{\lambda , c}(4)$. Therefore we have $l_\Omega (p_i)\le
s_{\lambda , c}(4)+s_{\lambda , c}(n-1)$. Repeating these
arguments for all $p_i\in I$, we get (\ref{ge5}).
\end{proof}

\begin{cor}
$s_{\lambda , c}(n)$ is finite for any $n$.
\end{cor}

To prove the main result of this section we need the following
auxiliary lemma. Although it is probably known, we did not find any
precise reference in the literature.

\begin{lem}\label{linfunct}
Let $f\colon \mathbb N\to \mathbb N$. Suppose that there exist
constants $C, N>0$, and $\alpha \in (0,1)$ such that for any $n\in
\mathbb N$, $n>N$, there are $n_1, \ldots , n_k\in \mathbb N$
satisfying the following conditions:

a) $k\le C \ln n$;

b) $f(n)\le \sum\limits_{i=1}^k f(n_i)$;

c) $n\le \sum\limits_{i=1}^k n_i\le n+C \ln n$;

d) $n_i\le \alpha n$  for any $i=1, \ldots , k$.

\noindent Then $f(n)$ is bounded by a linear function from above.
\end{lem}

\begin{proof}
Let $\e (n)=\ln \frac{C^2\ln ^2n}{1-\alpha } $ and let $N_0>N$ be
a constant such that
\begin{equation}\label{lin11}
2\e (n)\le \ln n
\end{equation}
and
\begin{equation}\label{lin12}
n-2C \ln n>0,
\end{equation}
for all $n\ge N_0$. Further let $N_1>N_0$ be a constant such that
\begin{equation}\label{lin13}
\frac{n(1-\alpha)}{C\ln n}\ge N_0
\end{equation}
for all $n\ge N_1$. The inequality (\ref{lin12}) allows us to
chose a positive constant $D$ such that
\begin{equation}\label{lin21}
f(n)\le Dn
\end{equation}
for all $n\le N_0$ and
\begin{equation}\label{lin22}
f(n)\le D(n-2C\ln n)
\end{equation}
for all $N_0<n\le N_1$. To prove the lemma it suffices to show that
(\ref{lin22}) holds for all $n\ge N_1$. We proceed by induction on
$n$.

Suppose that $n> N_1$. According to a) and c), there exists
$i_1\in \{ 1, \ldots , k\} $ such that $n_{i_1}\ge n/(C\ln n)$.
Furthermore, by d) we have $\sum\limits_{i\ne i_1} n_i\ge
n(1-\alpha)$. Hence there is $i_2\ne i_1$ such that $n_{i_2}\ge
n(1-\alpha )/(C\ln n)$. Note that $n_{i_j}\ge N_0$ for $j=1,2$ by
(\ref{lin13}). Obviously,
\begin{equation}\label{lin3}
\ln n_{i_1}+\ln n_{i_2} =\ln (n_{i_1}n_{i_2})\ge \ln
\frac{n^2(1-\alpha )}{C^2\ln ^2 n} = 2\ln n -\e (n).
\end{equation}

Let $\{ 1, \ldots , k\} =J_1\sqcup J_2,$ where $J_1$ consists of
all $j\in \{ 1, \ldots , k\} $ such that $n_j\ge N_0$. Applying
subsequently b), the inductive assumption together with
(\ref{lin21}) and (\ref{lin22}), c), (\ref{lin3}), and
(\ref{lin11}) we obtain
$$
\begin{array}{rl}
f(n)\le & \sum\limits_{i=1}^k f(n_i)\le \sum\limits_{i\in J_1}
D(n_i-2C\ln n_i) + \sum\limits_{i\in J_2} Dn_i \le
D\sum\limits_{i=1}^k n_i - 2CD (\ln n_{i_1}+\ln n_{i_2})\le\\ & \\
& D(n+C\ln n) -2CD(2\ln n -\e (n)) \le  D(n-2C\ln n) + DC(2\e (n)-\ln n)\le \\ & \\
& D(n-2C\ln n).
\end{array}
$$
\end{proof}

\begin{figure}
\hspace{-3mm} \includegraphics[width=160mm]{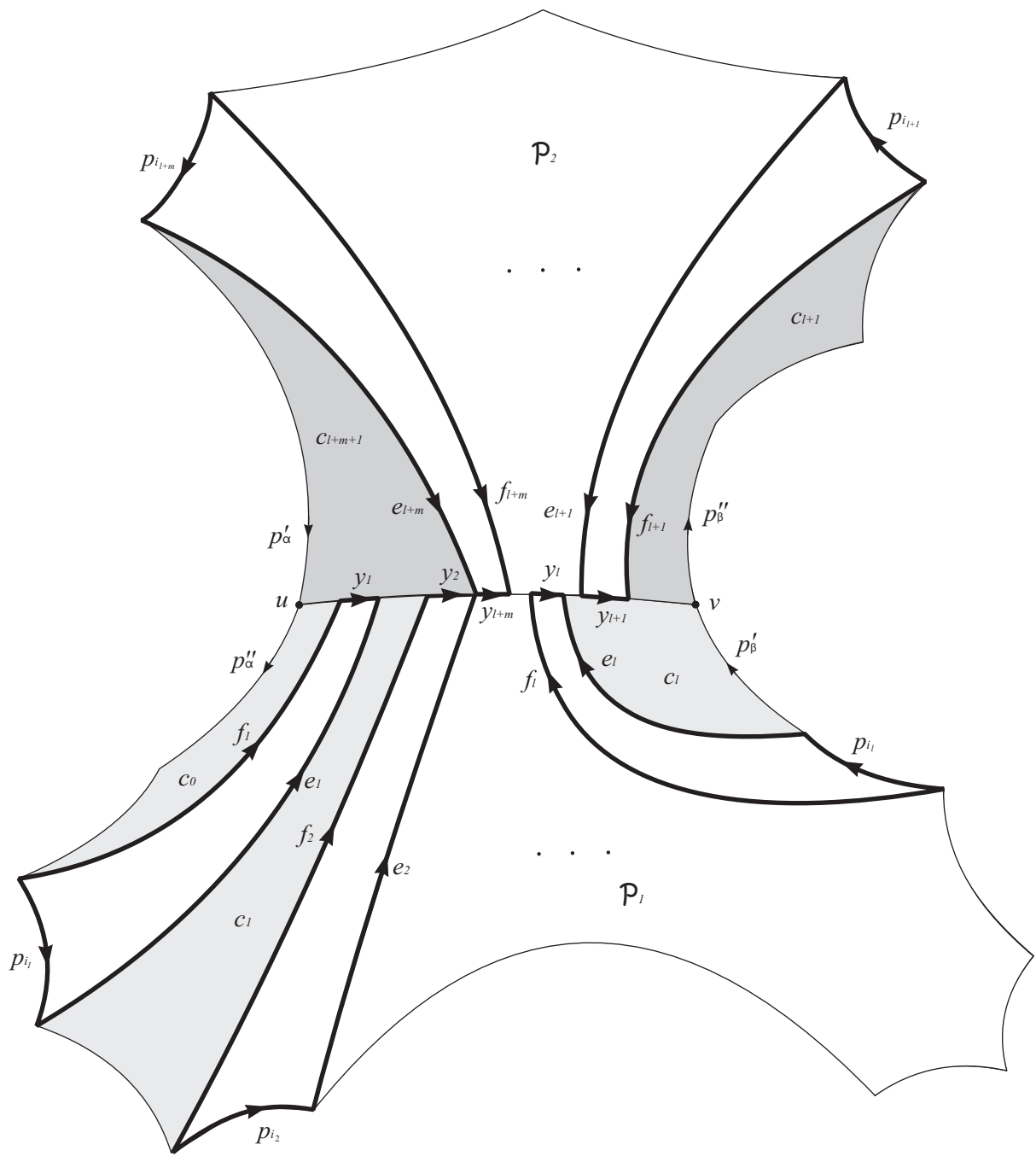}\\

\vspace{3mm}  \caption{Decomposition of the $n$--gon in the proof of
Proposition \ref{s(n)}} \label{fig5}
\end{figure}

Now we are ready to prove the main result of this section.

\begin{proof}[Proof of Proposition \ref{s(n)}]
We are going to show that for any fixed $\lambda \ge 1$, $c\ge 0$,
the function $s_{\lambda , c}(n)$ satisfies the assumptions of
Lemma \ref{linfunct}. Let $(\P , I)\in \mathcal Q_{\lambda ,
c}(n)$, where $\P =p_1\ldots p_n$. As in the proof of Lemma
\ref{l1}, we may assume that every $p_i\in I$ consists of a single
edge. We also assume $n\ge N$, where the constant $N$ is big
enough. The exact value of $N$ will be specified later.

Let $u,v$ be the points on $\P $ provided by Lemma \ref{P}.
Without loss of generality we may assume that $u,v$ are vertices
of $\G $. Further let $t$ denote a geodesic paths in $\G $ such
that $t_-=u$, $t_+=v$. According to Lemma \ref{P},
\begin{equation}\label{lr}
l(t)\le 2\delta (2+\log _2 n)+\theta ,
\end{equation}
where $\theta $ is a constant depending only on $\delta $,
$\lambda $, and $c$, and $t$ divides $\P $ into an $m_1$--gon $\P
_1$ and $m_2$--gon $\P _2$ such that
\begin{equation}\label{ni}
m_i\le 3n/4 +2<n
\end{equation}
for $i=1,2$. To be precise we assume that $u\in p_\alpha $, $v\in
p_\beta $, and $p_\alpha =p_\alpha ^\prime p_\alpha
^{\prime\prime}$,  $p_\beta =p_\beta ^\prime p_\beta
^{\prime\prime}$, where $(p_\alpha ^\prime)_+=(p_\alpha
^{\prime\prime})_-=u$, $(p_\beta ^\prime)_+=(p_\beta
^{\prime\prime})_-=v$. Then $$ \P _1 = p_\alpha ^{\prime\prime }
p_{\alpha +1} \ldots p_{\beta -1}p_\beta ^\prime t^{-1}$$ and $$
\P _2 = p_\beta ^{\prime\prime }p_{\beta +1}\ldots p_{\alpha
-1}p_\alpha ^\prime t.$$ (Here and below the indices are taken
modulo $n$.) Since each $p_i\in I$ consists of a single edge, one
of the paths $p_\alpha ^\prime $, $p_\alpha ^{\prime\prime }$
(respectively $p_\beta ^\prime $, $p_\beta ^{\prime\prime }$) is
trivial whenever $p_\alpha \in I$ (respectively $p_\beta \in I$).
Hence the set $I$ is naturally divided into two disjoint parts
$I=I_1\sqcup I_2$, where $I_i$ is a subset of $I$ consisting of
sides of $\P _i$, $i=1,2$.

Let us consider the polygon $\P _1$ and construct cycles $c_0,
\ldots , c_l$ in $\G $ as follows. If each $p_i\in I_1$ is isolated
in $\P _1$, we set $l=0$ and $c_0=\P _1$. Further suppose this is
not so. Let $p_{i_1}\in I_1$, be the first component (say, an
$H_{\lambda _1}$--component) in the sequence $p_\alpha , p_{\alpha
+1}, \ldots $ such that $p_{i_1}$ is not isolated in $\P _1$. As
$p_{i_1}$ is isolated in $\P $, this means that $p_{i_1}$ is
connected to an $H_{\lambda _1}$--component $y_1$ of $t$. Let $f_1$
(respectively $e_1$) be an edge in $\G $ labelled by an element of
$H_{\lambda _1}\setminus\{ 1\}$ (or a trivial path) such that
$(f_1)-=(p_{i_1})_-$, $(f_1)_+=(y_1)_-$ (respectively
$(e_1)_-=(p_{i_1})_+$, $(e_1)_+=(y_1)_+$). We set $$
c_0=p_\alpha^{\prime\prime }p_{\alpha +1}\ldots p_{i_1-1}f_1 [
(y_1)_-,u],$$ where $[(y_1)_-,u]$ is the segment of $t^{-1}$ (see
Fig. \ref{fig5}).

Now we proceed by induction. Suppose that the cycle $c_{k-1}$ and
the corresponding paths $f_{k-1},e_{k-1},y_{k-1}, p_{i_{k-1}}$ have
already been constructed. If the sequence $p_{i_{k-1} +1},
p_{i_{k-1} +2}, \ldots $ contains no component $p_i\in I_1$ that is
not isolated in $\P _1$, we set $l=k$,
$$ c_k=e_{k-1}^{-1}p_{i_{k-1}+1}\ldots p_{\beta -1}p_\beta ^\prime
[v,(y_{k-1})_+], $$ where $[v,(y_{k-1})_+]$ is the segment of
$t^{-1}$, and finish the procedure. Otherwise we continue as
follows. We denote by $p_{i_k}$ the first component in the
sequence $p_{i_{k-1} +1}, p_{i_{k-1} +2}, \ldots $ such that
$p_{i_k}\in I_1$  and $p_{i_k}$ is connected to some component
$y_k$ of $t$. Then we construct $f_k$, $e_k$ as above and set $$
c_k=e_{k-1}^{-1}p_{i_{k-1}+1}\ldots
p_{i_k-1}f_k[(y_k)_-,(y_{k-1})_+].$$ Observe that each path
$p_i\in I_1$ is either included in the set $J_1=\{ p_{i_1},
\ldots, p_{i_l}\} $ or is an isolated component of some $c_j$.
Indeed a paths $p_i\in I_1\setminus J_1$ can not be connected to a
component of $t$ according to our choice of $p_{i_1}, \ldots,
p_{i_l}$. Moreover $p_i\in I_1\setminus J_1$ can not be connected
to some $f_j$ or $e_j$ since otherwise $p_i$ is connected to
$p_{i_j}$ that contradicts the assumption that sides from the set
$I$ are isolated components in $\P $.

By repeating the 'mirror copy' of this algorithm for $\P _2$, we
construct cycles $c_{l+1}, \ldots , c_{l+m+1}$, $m\ge 0$, the set
of components $J_2=\{ p_{i_{l+1}},\ldots , p_{i_{l+m}}\}\subseteq
I_2$, components $y_{l+1}, \ldots , y_{l+m}$ of $t$, and edges (or
trivial paths) $f_{l+1}, e_{l+1}, \ldots,  f_{l+m}, e_{l+m}$ in
$\G $ such that $f_j$ (respectively $e_j$) goes from $(p_{i_j})_-$
to $(y_j)_+$ (respectively from $(p_{i_j})_+$ to $(y_j)_-$)  (see
Fig. \ref{fig5}) and each path $p_i\in I_2$ is either included in
the set $J_2$ or is an isolated component of $c_j$ for a certain
$j\in \{ l+1, \ldots, l+m+1\} $.

Each of the cycles $c_j$, $0\le j\le l+m+1$, can be regarded as a
geodesic $n_j$--gon whose set of sides consists of paths of the
following five types (up to orientation):

\begin{enumerate}
\item[(1)] Components from the set $I\setminus (J_1\cup J_2)$.

\item[(2)] Sides of $\P _1$ and $\P _2$ that do not belong to the
set $I$.

\item[(3)] Paths $f_j$ and $e_j$, $1\le j\le l+m$.

\item[(4)] Components $y_1 , \ldots , y_{l+m}$ of $t$.

\item[(5)] Maximal subpaths of $t$ lying 'between' $y_1, \ldots ,
y_{l+m}$, i.e. those maximal subpaths of $t$ that have no common
edges with $y_1, \ldots , y_{l+m}$.
\end{enumerate}

It is straightforward to check that for a given $0\le j\le l+m+1$,
all sides of $c_j$ of type (1), (3), and (4) are isolated
components of $c_j$. Indeed we have already explained that sides
of type (1) are isolated in $c_j$. Further, if $f_j$ or $e_j$ is
connected to $f_k$, $e_k$, or $y_k$ for $k\ne j$, then $p_{i_j}$
is connected to $p_{i_k}$ and we get a contradiction. For the same
reason $f_j$ or $e_j$ can not be connected to a component of a
side of type (2). If $f_j$ or $e_j$ is connected to a component
$x$ of a side of type (5), i.e., to a component of $t$, then $y_j$
is connected to $x$. This contradicts the assumption that $t$ is
geodesic. Finally $y_j$ can not be connected to a component of a
side of type (2) since otherwise $p_{i_j}$ is not isolated in $\P
$, and $y_j$ can not be connected to another component of $t$ as
notified in the previous sentence.

Observe that (\ref{lr}) and (\ref{ni}) imply the following
estimate of the number of sides of $c_j$:
$$
n_j\le \max\{ m_1, m_2\} + l(t)\le 3n/4 +2 +2\delta (\log _2 n +2)
+\theta .
$$
Assume that $N$ is a constant such that $3n/4 +2 +2\delta (\log _2
n +2)+\theta \le 4n/5 $ for all $n\ge N$. Then for any $n\ge N$,
we can apply the inductive assumption for the set of components of
type (1), (3), and (4) in each of the polygons $c_0, \ldots ,
c_{l+m+1}$. This yields
$$
\sum\limits_{p_i\in I} l_\Omega (p_i)\le \sum\limits_{p_i\in
I\setminus (J_1\cup J_2)} l_\Omega (p_i) +\sum\limits_{j=1}^{l+m}
\big( l_\Omega (y_j) + l_\Omega (e_j) +l_\Omega (f_j)\big) \le
\sum\limits_{j=0}^{l+m+1} s_{\lambda , c}(n_j)
$$
Further there is a constant $C>0$ such that
$$
\sum\limits_{j=0}^{m+l+1} n_j \le n + 6l(t)\le n +12\delta (\log
_2 n +2) +6\theta \le n+C\log _2 n
$$
and
$$
m+l+2\le 2l(t)+2\le C\log _2 n.
$$
Therefore, for any $n\ge N$, the function $s_{\lambda , c}(n)$
satisfies the assumptions of Lemma \ref{linfunct} for $k=m+l+2$
and $\alpha =4/5$. Thus $s(n, \lambda , c)$ is bounded by a linear
function from above.
\end{proof}

%%%%%%%%%%%%%%%%%%%%%%%%%%%%%%%%%%%%%%%%%%%%%%%%%%%%%%%%%%%%%%%

\section{Diagram surgery}

%%%%%%%%%%%%%%%%%%%%%%%%%%%%%%%%%%%%%%%%%%%%%%%%%%%%%%%%%%%%%%%

All conventions and notation from the beginning of the previous
section remain valid here. Together with the relative presentation
(\ref{G}) of $G$ with respect to $\Hl $ we also consider the
corresponding non-relative presentation (\ref{Gfull}). Given a
collection $\N =\{ N_\lambda \} _{\lambda \in \Lambda }$, where
$N_\lambda $ is a normal subgroup of $H_\lambda $, we denote by $N$
the normal closure of $\mbox{$\bigcup_{\lambda \in \Lambda }$}
N_\lambda $ in $G$. Recall that $G(\N )=G/N $. We fix the following
presentation for $G(\N )$
\begin{equation}\label{GNfull}
G(\N ) =\langle X\cup\mathcal H\; |\;\mathcal S\cup \mathcal Q\cup
\mathcal R \rangle ,
\end{equation}
where $\mathcal Q=\bigcup_{\lambda \in \Lambda } \mathcal
Q_\lambda $ and $\mathcal Q_\lambda $ consists of all words (not
necessary reduced) in the alphabet $H_\lambda \setminus \{ 1\} $
representing elements of $N_\lambda $ in $G$.

In this section we consider van Kampen diagrams over (\ref{Gfull})
of a certain type. More precisely, we denote by $\mathcal D$ the
set of all diagrams $\Delta $ over (\ref{Gfull}) such that:

(D1) Topologically $\Delta $ is a disc with $k\ge 0$ holes. More
precisely, the boundary of $\Delta $ is decomposed as $\partial
\Delta =\partial_{ext} \Delta \sqcup \partial_{int} \Delta $,
where $\partial _{ext}\Delta $ is the boundary of the disc and
$\partial _{int}\Delta $ consists of disjoint cycles ({\it
components}) $c_1, \ldots c_k$ that bound the holes.

(D2) For any $i=1, \ldots , k$, the label $\phi (c_i)$ is a word
in the alphabet $H_\lambda\setminus \{ 1\} $ for some $\lambda \in
\Lambda $ and this word represents an element of $N_\lambda $ in
$G$.

The following lemma relates diagrams of the described type to the
group $G(\N )$.

\begin{lem}\label{cutting}
A word $W$ in $X\cup \mathcal H$ represents $1$ in $G(\N )$ if and
only if there is a diagram $\Delta \in \mathcal D$ such that $\phi
(\partial _{ext} \Delta )\equiv W$.
\end{lem}

\begin{proof}
Suppose that $\Sigma $ is a disc van Kampen diagram over
(\ref{GNfull}). Then by cutting off all essential cells labeled by
words from $\mathcal Q$ ({\it $\mathcal Q$--cells}) and passing to a
$0$--refinement if necessary we obtain a van Kampen diagram $\Delta
\in \mathcal D$ with $\phi (\partial_{ext} \Delta )\equiv \phi
(\partial \Sigma )$. Conversely, each $\Delta \in D$ may be
transformed into a disk diagram over (\ref{GNfull}) by attaching
$\mathcal Q$--cells to all components of $\partial _{int} \Delta $.
\end{proof}

In what follows we also assume the diagrams from $\mathcal D$ to
be endowed with an additional structure.

(D3) Each diagram $\Delta \in \mathcal D$ is equipped with a {\it
cut system} that is a collection of disjoint paths ({\it cuts})
$T=\{ t_1, \ldots , t_k\} $ without self--intersections in $\Delta $
such that $(t_i)_+, (t_i)_-$ belong to $\partial \Delta $, and after
cutting $\Delta $ along $t_i$ for all $i=1, \ldots , k$ we get a
connected simply connected diagram $\widetilde{\Delta }$.

By $\kappa\colon \widetilde{\Delta }\to\Delta $ we denote the
natural map that 'sews' the cuts. We also fix an arbitrary point $O$
in $\widetilde{\Delta }$ and denote by $\mu $ the map provided by
Lemma \ref{map}.

\begin{lem}\label{trans}
Suppose that $\Delta \in \mathcal D$. Let $a,b$ be two vertices on
$\partial \Delta $, $\tilde a, \tilde b$ some vertices on
$\partial \widetilde{\Delta }$ such that $\kappa (\tilde a)=a$,
$\kappa (\tilde b)= b$. Then for any paths $r$ in $\G $ such that
$r_-=\mu (\tilde a)$, $r_+=\mu (\tilde b)$, there is a diagram
$\Delta _1\in \mathcal D$ endowed with a cut system $T_1$ such
that the following conditions hold:
\begin{enumerate}
\item $\Delta _1$ has the same boundary and the same cut system as $\Delta $.
By this we mean the following. Let $\Gamma _1$ (respectively $\Gamma
$) be the subgraph of the $1$-skeleton of $\Delta_1$ (respectively
of the $1$-skeleton of $\Delta $) consisting of $\partial \Delta _1$
(respectively $\partial \Delta $) and all cuts from $T_1$
(respectively $T$). Then there is a graph isomorphism $\Gamma _1\to
\Gamma $ that preserves labels and orientation and maps cuts of
$\Delta _1$ to cuts of $\Delta $ and $\partial _{ext} \Delta _1$ to
$\partial _{ext} \Delta $.

\item There is a paths $q$ in $\Delta _1$ without
self--intersections such that $q_-=a$, $q_+=b$, $q$ has no common
vertices with cuts $t\in T_1$ except for possibly $a$,$b$, and $\phi
(q)\equiv \phi (r)$.
\end{enumerate}
\end{lem}

\begin{figure}
\vspace{1mm}
\hspace{1.4cm}\includegraphics[]{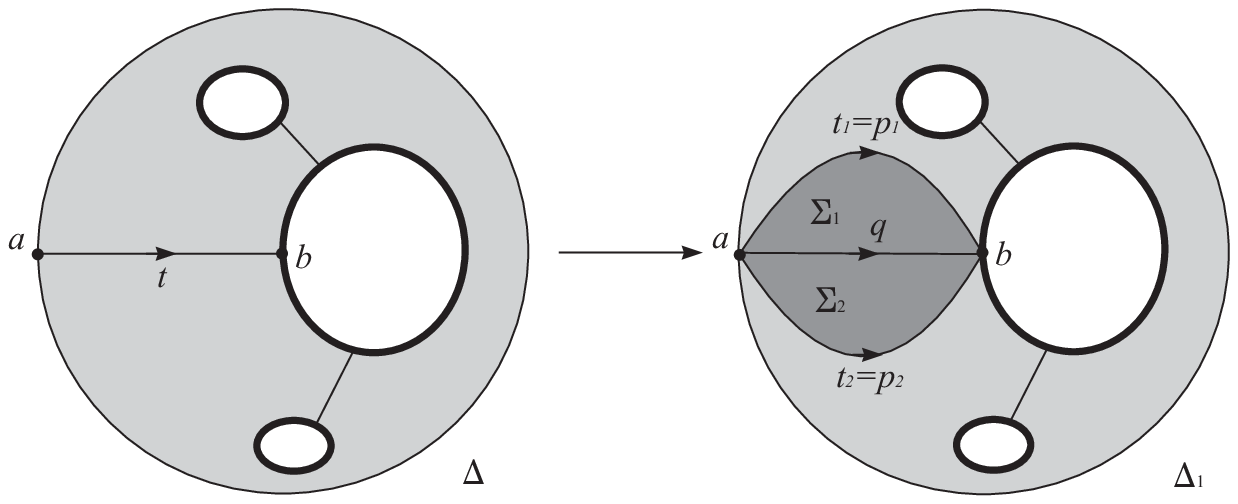}\\
\vspace{-2mm}
  \caption{}\label{fig6}
\end{figure}

\begin{proof}
Let us fix an arbitrary path $\tilde t$ in $\widetilde{\Delta }$
without self--intersections that connects $\tilde a$ to $\tilde b$
and intersects $\partial \widetilde{\Delta }$ at the points $\tilde
a$ and $\tilde b$ only. The last condition can always be ensured by
passing to a $0$--refinement of $\Delta $ and the corresponding
$0$--refinement of $\widetilde{\Delta }$. Thus $t=\kappa (\tilde t)$
connects $a$ to $b$ in $\Delta $ and has no common points with cuts
$t\in T$ except for possibly $a$,$b$. Note that
$$\phi (t)\equiv \phi (\tilde t)\equiv \phi(\mu (\tilde t))$$
as both $\kappa $, $\mu $ preserve labels and orientation.

Since  $\mu (\tilde t)$ connects $\mu (\tilde a)$ to $\mu (\tilde
b)$ in $\G$, $\phi (\mu (\tilde t))$ represents the same element
of $G$ as $\phi (r)$. Hence there exists a disk diagram $\Sigma
_1$ over (\ref{Gfull}) such that $\partial \Sigma _1=p_1q^{-1}$,
where $\phi (p_1)\equiv \phi (t)$ and $\phi (q)\equiv \phi (r)$.
Let $\Sigma _2$ denote its mirror copy. We glue $\Sigma _1$ and
$\Sigma _2$ together by attaching $q$ to its mirror copy. Thus we
get a new diagram $\Sigma $ with boundary $p_1p_2^{-1}$, where
$\phi (p_1)\equiv\phi (p_2)\equiv \phi (t)$. The path in $\Sigma $
corresponding to $q$ in $\Sigma _1$ and its mirror copy in $\Sigma
_2$ is also denoted by $q$.

We now perform the following surgery on the diagram $\Delta $.
First we cut $\Delta $ along $t$ and denote the new diagram by
$\Delta _0$. Let $t_1$ and $t_2$ be the two copies of the path $t$
in $\Delta _0$. Then we glue $\Delta _0$ and $\Sigma $ by
attaching $t_1$ to $p_1$ and $t_2$ to $p_2$ (Fig. \ref{fig6}) and
get a new diagram $\Delta _1$. This surgery does not affect cuts
of $\Delta $ as $t$ had no common points with cuts from $T$ except
for possibly $a$ and $b$. Thus the system of cuts in $\Delta _1$
is inherited from $\Delta $ and $\Delta _1$ satisfies all
requirements of the lemma.
\end{proof}

\begin{defn}
By an {\it $H_\lambda $--path} in $\Delta \in \mathcal D$ or in
$\widetilde{\Delta }$ we mean any paths whose label is a
(nontrivial) word in $H_\lambda \setminus \{ 1\} $. We say that two
such paths $p$ and $q$ in $\Delta \in \mathcal D$ are {\it
connected} if they are $H_\lambda $--paths for the same $\lambda \in
\Lambda $ and there are $H_\lambda $--paths $a$, $b$ in
$\widetilde{\Delta }$ such that $\kappa (a)$ is a subpaths of $p$,
$\kappa (b)$ is a subpaths of $q$, and $\mu (a)$, $\mu (b)$ are
connected in $\G $, i.e., there is a path in $\G $ that connects a
vertex of $\mu (a)$ to a vertex of $\mu (b)$ and is labelled by a
word in $H_\lambda \setminus \{ 1\} $. We stress that the equalities
$\kappa (a)=p$ and $\kappa (b)=q$ are not required. Thus the
definition makes sense even if the paths $p$ and $q$ are cut by the
cuts of $\Delta $ into several pieces.
\end{defn}

\begin{defn}\label{typeD}
We also define the {\it type} of a diagram $\Delta \in
\mathcal D$ by the formula
$$
\tau (\Delta )=\left( k, \sum\limits_{i=1}^k l(t_i)\right) ,
$$
where $k$ is the number of holes in $\Delta $. We fix the standard
order on the set of all types by assuming $(m,n)\le (m_1, n_1)$ is
either $m< m_1$ or $m=m_1 $ and $n\le n_1$.
\end{defn}

For a word $W$ in the alphabet $X\cup\mathcal H$, let $\mathcal
D(W)$ denote the set of all diagrams $\Delta \in \mathcal D$ such
that $\phi(\partial _{ext} \Delta )\equiv W$. In the proposition
below we say that a word $W$ in $X\cup \mathcal H$ is {\it
geodesic} if any (or, equivalently, some) path in $\G $ labelled
by $W$ is geodesic.

\begin{prop}\label{DS}
Suppose that $W$ is a word in $X\cup\mathcal H$ representing $1$
in $G(\N )$, $\Delta $ is a diagram of minimal type in $\mathcal
D(W)$, $T$ is the cut system in $\Delta $, and $c$ is a component
of $\partial _{int}\Delta $. Then:
\begin{enumerate}
\item For each cut $t\in T$, the word $\phi (t)$ is geodesic.

\item The label of $c$ represents a nontrivial element in $G$.

\item The path $c$ can not be connected to an $H_\lambda
$--subpath of a cut.

\item The path $c$ can not be connected to another component of
$\partial _{int} \Delta $
\end{enumerate}
\end{prop}

\begin{figure}
 \hspace{3mm}
 \includegraphics[]{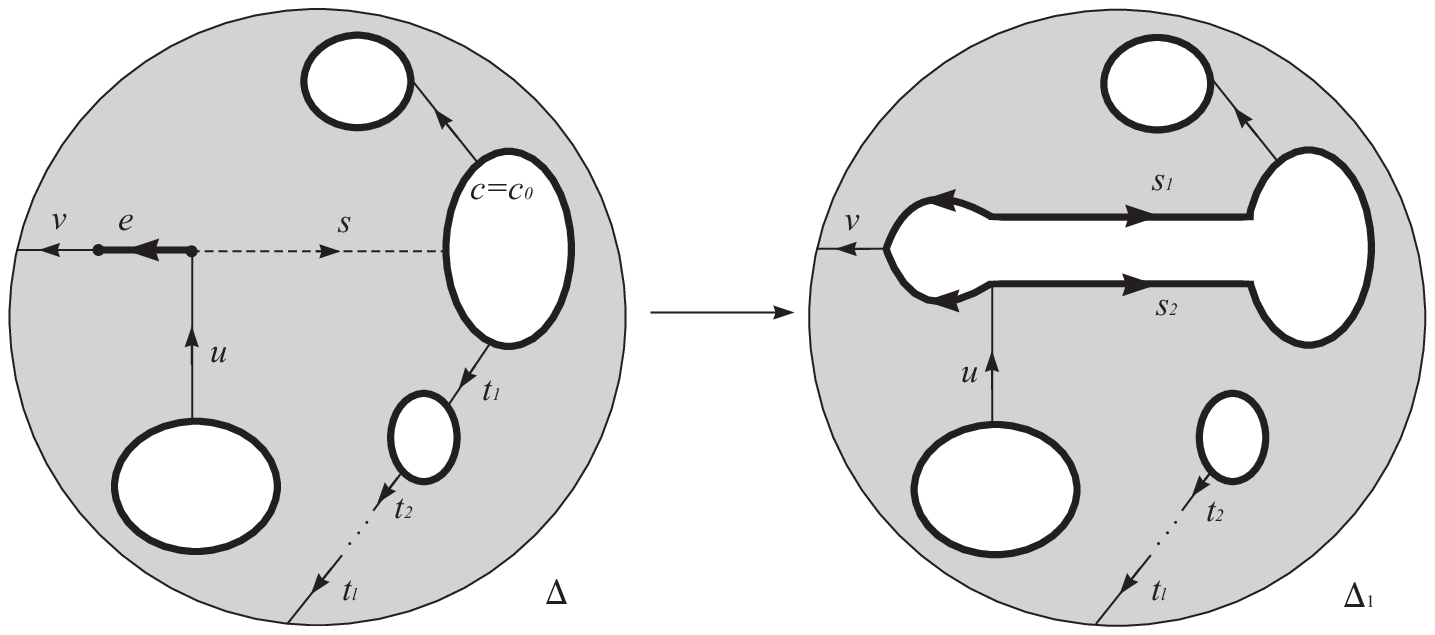}

\begin{center} a) \end{center}
\vspace{5mm}

\hspace{3mm}
\includegraphics[]{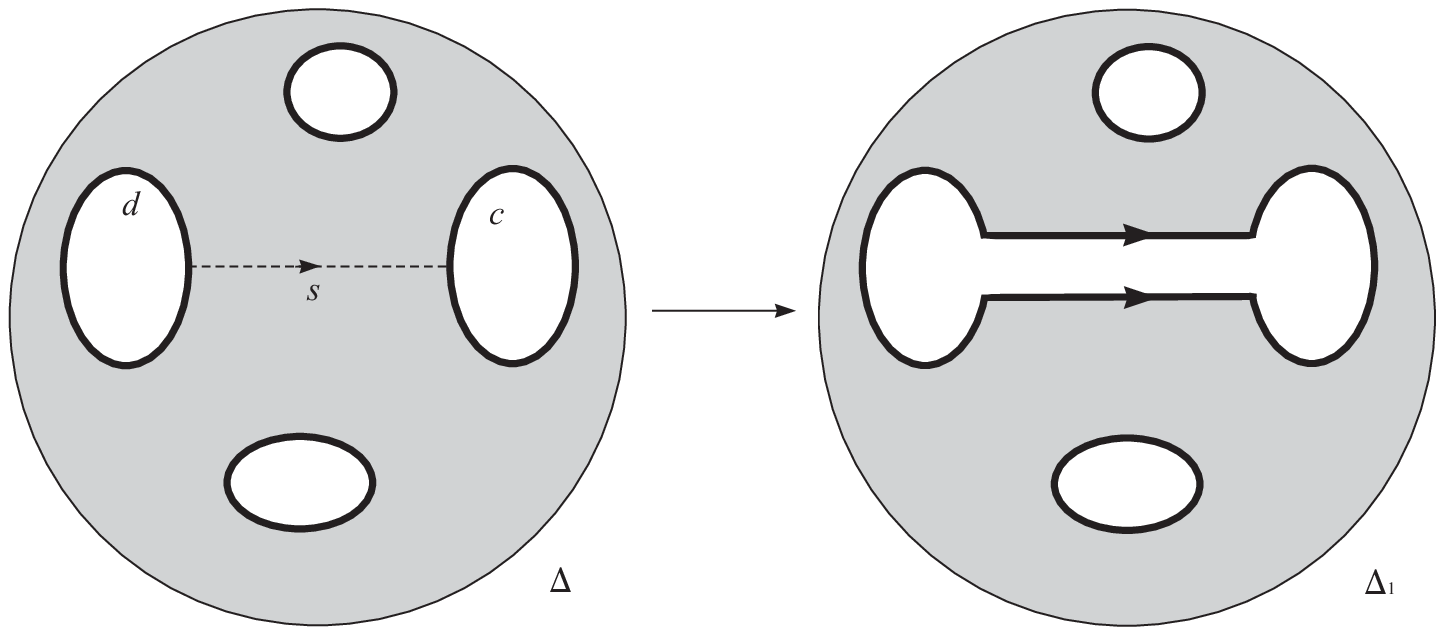}\\

\begin{center} b) \end{center}

  \caption{} \label{fig7}
\end{figure}

\begin{proof}
Assume that for a certain path $t\in T$, $\phi (t)$ is not geodesic.
Let $\tilde a$, $\tilde b$ be vertices in $\widetilde \Delta $ such
that $\kappa (\tilde a)=t_-$, $\kappa (\tilde b)=t_+$. Let also $r$
be a geodesic paths in $\G $ that connects $\mu (\tilde a)$ to $\mu
(\tilde b)$. Applying Lemma \ref{trans}, we may assume that there is
a path $q$ in $\Delta $ such that $q_-=t_-$, $q_+=t_+$, and $\phi
(q)\equiv \phi (r)$, i.e., $\phi (q)$ is geodesic. In particular,
$l(q)<l(t)$. Now replacing $t$ with $q$ in the cut system we reduce
the type of the diagram. This contradicts the choice of $\Delta $.

The second assertion is obvious. Indeed if $\phi (c)$ represents $1$
in $G$, there is a disk diagram $\Pi $ over (\ref{Gfull}) with
boundary label $\phi (\partial \Pi)\equiv \phi (c)$. Attaching $\Pi
$ to $c$ does not affect $\partial _{ext}\Delta $ and reduces the
number of holes in the diagram. This contradicts the minimality of
$\tau (\Delta )$ again.

Further assume that $c$ is connected to an $H_\lambda $--subpath $e$
of some $r\in T$. Then $c$ is an $H_\lambda $--path for the same
$\lambda \in \Lambda$. Let $r=uev$. Cutting $\Delta $ along $e$ (to
convert $e$ into a boundary component), applying Lemma \ref{trans},
and gluing the copies of $e$ back, we may assume that there is a
path $s$ without self--intersections in $\Delta $ such that
$s_-=e_-$, $s_+\in c$, and $\phi (s)$ is a word in $H_\lambda
\setminus \{ 1\} $. Moreover passing to a $0$--refinement, we may
assume that $s$ has no common vertices with the boundary of the
diagram, paths from $T\setminus \{ r\} $, $u$, and $v$ except for
$s_-$ and  $s_+$. Now we cut $\Delta $ along $s$ and $e$. Let $s_1$,
$s_2$ be the copies of $s$ in the obtained diagram $\Delta _1$. The
boundary component of $\Delta _1$ obtained from $c$ and $e$ has
label $\phi(c)\phi (s)^{-1}\phi(e)\phi(e)^{-1}\phi{s}$ that is a
word in $H_\lambda \setminus\{ 1\} $ representing an element of
$N_\lambda $ in $G$. Note also that our surgery does not affect cuts
of $\Delta $ except for $r$. Thus the system of cuts $T_1$ in
$\Delta _1$ may obtained from $T$ as follows. Since $\widetilde
{\Delta } $ is connected and simply connected, there is a unique
sequence
$$
c=c_0,\; t_1,\; c_1,\; \ldots ,\; t_l,\; c_l=\partial _{ext} {\Delta
},
$$
where $c_0, \ldots , c_l$ are (distinct) components of $\partial
\Delta $, $t_i\in T$, and (up to orientation) $t_i$ connects
$c_{i-1} $ to $c_i$, $i=1, \ldots , l$ (Fig. \ref{fig7}a). We set
$T_1=(T\setminus \{ r, t_1\} )\cup \{ u, v\} $. Thus $\Delta _1\in
\mathcal D(W)$ and $\tau (\Delta _1)<\tau (\Delta )$. Indeed
$\Delta _1$ and $\Delta $ have the same number of holes and
$\sum\limits_{t\in T_1} l(t)\le \sum\limits_{t\in T_1} l(t)-1$.
This contradicts the choice of $\Delta $.

Finally suppose that $c$ is connected to another component $d$ of
$\partial _{int} \Delta $, $d\ne c$. To be definite, assume that $c$
and $d$ are labelled by words in $H_\lambda \setminus\{ 1\} $. Again
without loss of generality we may assume that there is a path $s$
without self--intersections in $\Delta $ such that $s_-\in d$,
$s_+\in c$, $\phi (s)$ is a word in $H_\lambda \setminus \{ 1\} $,
and $s$ has no common points with $\partial \Delta $ and paths from
$T$ except for $s_-$ and $s_+$. Let us cut $\Delta $ along $s$ and
denote by $\Delta _1$ the obtained diagram (Fig. \ref{fig7}b). This
transformation does not affect $\partial _{ext} \Delta $ and the
only changed internal boundary component has label $\phi
(c)\phi(s)^{-1}\phi (d)\phi (s)$, which is a word in $H_\lambda
\setminus\{ 1\} $. This word represents an element of $N_\lambda $
in $G$ as $N_\lambda\lhd H_\lambda $. We now fix an arbitrary system
of cuts in $\Delta _1$. Then $\Delta _1\in \mathcal D(W)$ and the
number of holes in $\Delta _1$ is smaller that the number of holes
in $\Delta $. We get a contradiction again.
\end{proof}

%%%%%%%%%%%%%%%%%%%%%%%%%%%%%%%%%%%%%%%%%%%%%%%%%%%%%%%%%

\section{Proofs of the main results}

%%%%%%%%%%%%%%%%%%%%%%%%%%%%%%%%%%%%%%%%%%%%%%%%%%%%%%%%%

Recall that $\Omega $ denotes the set provided by Lemma
\ref{Omega}. Let $D=D(2,0)$ be the constant from Proposition
\ref{s(n)}. We set
$$\mathcal F=\{ g\in \langle \Omega \rangle \; :\; |g|_\Omega \le
4D\} \setminus\{ 1\}. $$ Throughout the rest of the section we
assume that $N_\lambda \cap \mathcal F =\emptyset $ for all
$\lambda \in \Lambda $.

\begin{defn}\label{theta}
Given a word $W$ in the alphabet $X\cup \mathcal H$ representing
$1$ in $G(\N )$, we denote by $q(W)$ the minimal number of holes
among all diagrams from $\mathcal D(W)$. Further we define the
{\it type} of $W$ by the formula $\theta (W)=(q(W), \| W\| )$. The
set of types is endowed with the natural order (as in Definition
\ref{typeD}).
\end{defn}

The proof of Theorem \ref{CEP} is divided into a sequence of
lemmas. We begin with the first assertion of the theorem. Recall
that a word $W$ in $X\cup \mathcal H$ is called $(\lambda ,
c)$--quasi--geodesic (in $G$) for some $\lambda \ge 1$, $c\ge 0$,
if some (or, equivalently, any) path in $\G $ labelled by $W$ is
$(\lambda , c)$--quasi--geodesic. The following three results are
proved by common induction on $q(W)$.

\begin{lem}\label{CEP1}
Suppose that $W$ is a word in the alphabet $X\cup \mathcal H$
representing $1$ in $G(\N )$ and $\Delta $ is a diagram of minimal
type in $\mathcal D (W)$. Then:

\begin{enumerate}
\item Assume that for some $\lambda \in \Lambda $, $p$ and $q$ are
two connected $H_\lambda $--subpaths of the same component $c$ of
$\partial _{int}\Delta $, then there is an $H_\lambda $--component
$r$ of $\partial\widetilde{\Delta }$ such that $p$ and $q$ are
subpaths of $\kappa (r)$.

\item  If $W$ is $(2,0)$--quasi--geodesic and $q(W)>0$, then some
component of $\partial _{int} \Delta $ is connected to an
$H_\lambda $--subpath of $\partial _{ext}\Delta $ for some
$\lambda \in \Lambda $.

\item If $W$ is a word in the alphabet $H_\lambda \setminus \{ 1\}
$ for some $\lambda \in \Lambda $, then $W$ represents an element
of $N_\lambda $ in $G$.
\end{enumerate}
\end{lem}

\begin{proof}
For $q(W)=0$ the lemma is trivial. Assume that $q(W)>0$.

\begin{figure}
  \vspace{1mm}
 \hspace{13mm}\includegraphics[]{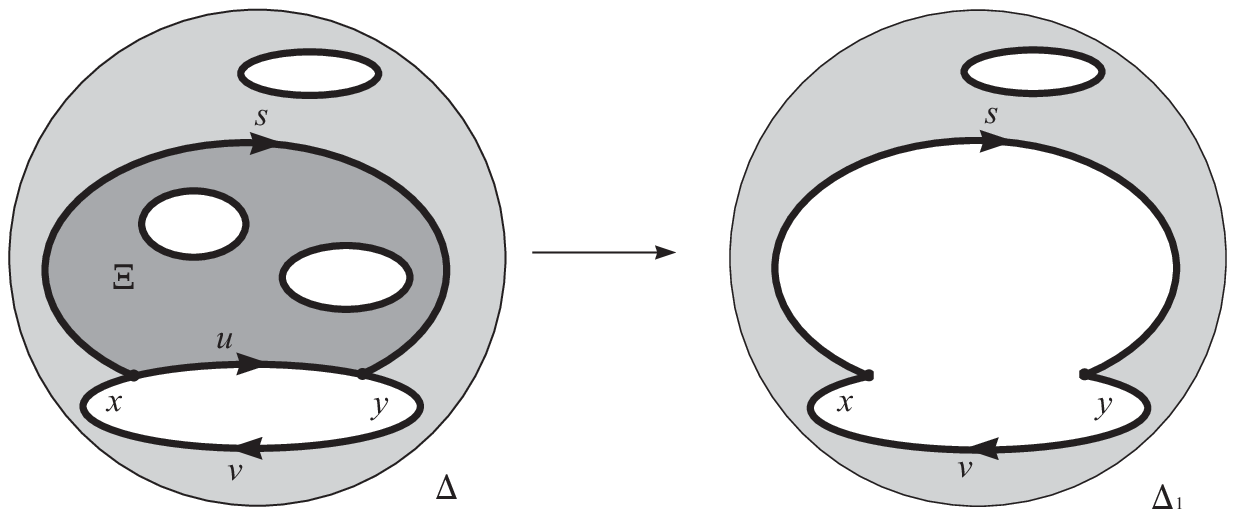}
  \caption{}\label{9}
\end{figure}

Let us prove the first assertion. Let $x$ (respectively $y$) be an
ending vertex of a certain essential edge of $p$ (respectively
$q$). Passing to a $0$--refinement of $\Delta $, we may assume
that $x$ and $y$ do not belong to any cut from the cut system $T$
of $\Delta $. Applying Lemma \ref{trans} we get a paths $s$ in
$\Delta $ connecting $x$ to $y$ such that $\phi (s)$ is a word in
the alphabet $H_\lambda \setminus \{ 1\}$ and $s$ does not
intersect any path from $T$. Let us denote by $\Xi $ the
subdiagram of $\Delta $ bounded by $s$ and the segment $u=[x,y]$
of $c^{\pm 1}$ such that $\Xi $ does not contain the hole bounded
by $c$ (Fig. \ref{9}).

Note that $V\equiv \phi (\partial \Xi )$ is a word in the alphabet
$H_\lambda \setminus \{ 1\}$ and $q(V)< q(W)$. By the third
assertion of our lemma, $V$ represents an element of $N_\lambda $.
Up to a cyclic shift, the label of the external boundary component
of the subdiagram $\Sigma =\Xi \cup c$ of $\Delta $ is a word in
$H_\lambda \setminus \{ 1\}$ representing the same element as
$\phi (c^{\pm 1})\phi(u)^{-1}V^{\pm 1}\phi (u)$ in $G$. As
$N_\lambda $ is normal in $H_\lambda $ and $\phi (u)$ represents
an element of $H_\lambda $ in $G$, $\phi (\partial _{ext} \Sigma
)$ represents an element of $N_\lambda $ in $G$. If $\Xi $
contains at least one hole, we replace $\Sigma$ with a single hole
bounded by $\partial _{ext}\Sigma $ (Fig. \ref{9}). This reduces
the number of holes in $\Delta $ and we get a contradiction.
Therefore $\Xi $ is simply connected. In particular, the path $u$
does not intersect any cut from $T$. This means that $p$ and $q$
are covered by the image of the same $H_\lambda $--component of
$\partial\widetilde{\Delta }$.

To prove the second assertion we suppose that for every component
$c_i$ of $\partial _{int} \Delta $, no $H_\lambda $--subpath of
$\partial _{ext}\Delta $ is connected to $c$. Then Proposition
\ref{DS} and the first assertion of our lemma imply that each
component $c_i$ of $\partial _{int}\Delta $ gives rise to
$H_\lambda $--components $a_{i1}, \ldots, a_{il}$ of $\partial
\widetilde{\Delta }$  for some $l=l(i)$ such that $\kappa
(a_{ij})\in c_i$, $j=1, \ldots , l$, and $\mu (a_{i1}),\ldots ,\mu
(a_{il})$ are isolated $H_\lambda $--components of the cycle
$\mathcal P=\mu (\partial \widetilde{\Delta })$ in $\G $.

Recall that for an element $g\in G$, the $\Omega $--lengths
$|g|_\Omega $ is defined to be the word lengths of $g$ with
respect to $\Omega $ if $g\in \langle \Omega \rangle $ and $\infty
$ otherwise. For each component $c_i$ of $\partial _{int} (\Delta
)$, we fix a vertex $o_i\in c_i$ such that $o_i=t_-$ or $o_i=t_+$
for some $t\in T$ and denote by $g_i$ the element represented by
$\phi (c_i)$ when we read this label starting from $o_i$. Clearly
\begin{equation}\label{gi}
|g_i| _{\Omega } \le \sum\limits_{j=1}^{l(i)} l_\Omega (\mu
(a_{ij})).
\end{equation}

The path $\mathcal P$ may be considered as an $n\le 4q(W)$--gon
whose sides (up to orientation) are of the following three types:
\begin{enumerate}
\item[(1)] sides corresponding to parts of $\partial _{ext} \Delta
$;

\item[(2)] sides corresponding to cuts in $\Delta $;

\item[(3)] components corresponding to $\partial _{int}\Delta $.
\end{enumerate}

The sides of $\mathcal P$ of type (1) are $(2,0)$--quasi--geodesic
in $\G $ as $W$ is $(2,0)$--quasi--geodesic. The sides of type (2)
are geodesic in $\G $ by the first assertion of Proposition
\ref{DS}. Hence we may apply Proposition \ref{s(n)} to the
$n$--gon $\mathcal P$, where the set of components $I$ consists of
sides of type (3). Taking into account (\ref{gi}), we obtain
$$
\sum\limits_{i=1}^{q(W)} |g_i|_\Omega \le \sum\limits_{p\in
I}l_{\Omega } (p) \le Dn\le 4D q(W),
$$
where $D=D(2,0)$ is provided by Proposition \ref{s(n)}.  Hence at
least one element $g_i$ satisfies $|g_i|_\Omega <4D $. According to
our choice of $\mathcal F$ and $\N $, we have $g_i=1$ in $G$.
However this contradicts the second assertion of Proposition
\ref{DS}.

To prove the last assertion we note that it suffices to deal with
the case when $W$ is geodesic as any element of $H_\lambda $ can
be represented by a single letter. Let $\Delta $ be a diagram of
minimal type in $\mathcal D (W)$. By the second assertion of the
lemma, some component $c$ of $\partial _{int}\Delta $ labelled by
a word in $H_\lambda \setminus \{ 1\} $ is connected to
$\partial_{ext}\Delta $. Applying Lemma \ref{trans} yields a path
$s$ in $\Delta $ connecting $\partial_{ext}\Delta $ to $c$ such
that $\phi (s) $ is a word in the alphabet $H_\lambda \setminus \{
1\}$. Let us cut $\Delta $ along $s$ and denote the new diagram by
$\Delta _1$. Obviously the word
$$
\phi (\Delta _1)\equiv \phi (s)\phi(c)\phi (s^{-1})\phi
(\Delta)
$$
is a word in the alphabet $H_\lambda \setminus \{ 1\}$ and $q(\phi
(\Delta _1))<q(W)$. By the inductive assumption, $\phi (\Delta
_1)$ represents an element of $N_\lambda $ in $G$. Since $\phi
(c)$ represents an element of $N_\lambda $ and $N_\lambda \lhd
H_\lambda $, the word $\phi (\Delta )$ also represents an element
of $N_\lambda $.
\end{proof}

The third assertion of Lemma \ref{CEP1} obviously implies the
first assertion of Theorem \ref{CEP}. Let us prove the second one.
For a word $W$ in the alphabet $X\cup \mathcal H$ representing $1$
in $G(\N )$, we set $$\AA (W)=\min\limits_{\Delta \in \mathcal
D(W)} N_\mathcal R (\Delta ).$$ It is easy to see that for any two
words $U$ and $V$ in $X\cup \mathcal H$ representing $1$ in $G(\N
)$, we have
\begin{equation}\label{AA}
\AA (UV)\le \AA (U) +\AA (V).
\end{equation}

\begin{lem}\label{CEP2}
For any word $W$ in $X\cup \mathcal H$ representing $1$ in $G(\N
)$, we have $\AA (W)\le 3C\| W\| $, where $C$ is the relative
isoperimetric constant of (\ref{G}).
\end{lem}

\begin{proof}
The proof is by induction on $\theta (W) $ (see Definition
\ref{theta}). If $q(W)=0$, then $W=1$ in $G$ and the required
estimate on $\AA (W)$ follows from the relative hyperbolicity of
$G$. We now assume that $q(W)>1$.

First suppose that the word $W$ is not $(2,0)$--quasi--geodesic in
$G$. That is, up to a cyclic shift $W\equiv W_1W_2$, where $W_1=U$
in $G$ and $\| U\| < \| W_1\| /2$. Note that $q(W_1U^{-1})=0$,
$q(UW_2)=q(W)$, and $\| UW_2\| \le \| W\| -\| W_1\| /2$. Hence
$\theta (UW_2)<\theta (W)$. Using the inductive assumption and
(\ref{AA}), we obtain
$$
\begin{array}{rl}
\AA (W)\le & \AA (W_1U^{-1})+\AA (UW_2)< \\ & \\ & \frac32C \|
W_1\| + 3C \left(\| W\| -\frac12 \| W_1\| \right) = 3C\| W\|.
\end{array}
$$

Now assume that $W$ is $(2,0)$--quasi--geodesic. Let $\Delta $ be
a diagram of minimal type in $\mathcal D(W)$. By the second
assertion of Lemma \ref{CEP1}, some component $c$ of $\partial
_{int} \Delta $ is connected to an $H_\lambda $--subpath $p$ of
$\partial _{ext} \Delta $ for some $\lambda \in \Lambda $.
According to Lemma \ref{trans}, we may assume that there is a path
$s$ in $\Delta $ connecting $c$ to $p_+$ such that $\phi (s)$ is a
word in the alphabet $H_\lambda \setminus \{ 1\} $. We cut $\Delta
$ along $s$ and denote by $\Delta _1$ the obtained diagram. Up to
cyclic shift, we have $W\equiv W_0\phi (p)$ and $$\phi (\Delta
_1)\equiv W_0\phi (p)\phi (s)^{-1}\phi (c)\phi (s).$$ Let $h$ be
the element of $H_\lambda $ represented by $\phi (p)\phi
(s)^{-1}\phi (c)\phi (s)$ in $G$. Observe that $q(W_0h)=q(\phi
(\Delta _1))<q(W)$. Further since $h^{-1}\phi (p)$ is a word in
$H_\lambda \setminus \{ 1\} $ representing $1$ in $G(\N )$, we
have $h^{-1}\phi (p)\in \mathcal Q$ and hence $\AA (h^{-1}\phi
(p))=0$. Applying the inductive assumption we obtain
$$
\begin{array}{rl}
\AA (W)= & \AA (W_0h)+\AA (h^{-1}\phi (p)) = \\ & \\ & \AA (W_0h)
\le 3C \| W_0h\| \le 3C \| W\| .
\end{array}
$$
\end{proof}

It is easy to see that the second assertion of Theorem \ref{CEP}
follows from Lemma \ref{CEP2}. Indeed, let $\e_1 \colon F(\N )\to
G(\N )$ be the natural homomorphism, where $F(\N )=F(X)\ast (\ast
_{\lambda \in \Lambda }H_\lambda /N_\lambda )$. Let $\e _0$ denote
the natural homomorphism $F\to F(\N )$, where $F$ is given by
(\ref{F}). The first assertion of Theorem \ref{CEP} implies that
$Ker\, \e _1=\langle \e _0 (\mathcal R)\rangle ^{F(\N )}$. Now let
$U$ be an element of $F(\N )$ such that $\e _1 (U)=1$, $W\in F$ a
preimage of $U$ such that $\| W\| =\| U\| $. Lemmas \ref{CEP2} and
\ref{cutting} imply that
\begin{equation}\label{WW1}
W=_F\prod\limits_{i=1}^{k} f_i^{-1}R_i^{\pm 1}f_i,
\end{equation}
where $f_i\in F$, $R_i\in \mathcal R\cup \mathcal Q$, and the
number of multiples corresponding to $R_i\in \mathcal R$ is at
most $3C\| W\| $. Applying $\e _0 $ to the both sides of
(\ref{WW1}) and taking into account that $\e_0 (f_i^{-1}R_if_i)=1$
in $F(\N )$ whenever $R_i\in \mathcal Q$, we obtain
$$
U=_{F(\N )} \prod\limits_{i=1}^{l} g_i^{-1}P_i^{\pm 1}g_i,
$$
where $g_i\in F(\N )$, $P_i\in \e _0 (\mathcal R)$, and $l\le 3C
\| W\| =3C\| U\| $. By definition this means that $G(\N )$ is
hyperbolic relative to $\{ H_\lambda /N_\lambda \} _{\lambda \in
\Lambda }$.

Let us prove the last assertion of the theorem. Since $S$ is
finite, without loss of generality we may assume that for any two
elements $s,t\in S$, we have $st^{-1}\in X$. Thus it suffices to
prove the following.

\begin{lem}
For any element $x\in X$, $x=1$ in $G(\N )$ implies $x=1$ in $G$.
\end{lem}

\begin{proof}
Suppose that $x=1$ in $G(\N )$ for some $x\in X$. Assume that
$x\ne 1$ in $G$. Then $q(x)>0$. Let $\Delta $ be a diagram of
minimal type in $\mathcal D (x)$. Since $x$ is a geodesic word in
$G$, some component of $\partial _{int}\Delta $ is connected to an
$H_\lambda $--subpath of $\partial _{ext}\Delta $ for some
$\lambda \in \Lambda $ by the second assertion of Lemma
\ref{CEP1}. However $\partial _{ext} \Delta $ contains no
$H_\lambda $--subpaths at all and we get a contradiction.
\end{proof}

Finally, to prove Corollary \ref{frh}, we need some results about
elementary subgroups of relatively hyperbolic groups obtained in
\cite{ESBG}. Recall that an element $g\in G$ is called {\it
hyperbolic} if it is not conjugate to an element of one of the
subgroups $H_\lambda $, $\lambda \in \Lambda $.

\begin{lem}\label{E(g)}
Let $g$ be a hyperbolic element of infinite order in $G$. Then
\begin{enumerate}
\item The element $g$ is contained in a unique maximal elementary
subgroup $E_G(g)$ of $G$.

\item The group $G$ is hyperbolic relative to the collection
$\Hl\cup \{ E_G(g)\} $.
\end{enumerate}
\end{lem}

The next result is also proved in \cite[Corollary 4.5]{ESBG}. (For
finitely generated group it can also be proved by using the action
of $G$ on its boundary, see \cite{Tuk}.)

\begin{lem}\label{hypel}
Suppose that all subgroups $\Hl $ are proper. Then $G$ contains a
hyperbolic element of infinite order.
\end{lem}

Recall that two elements $f,g\in G$ are called {\it commensurable}
if $f^k$ is conjugate to $g^l$ in $G$ for some $k,l\ne 0$. The lemma
below is a particular case of \cite[Theorem 1.4]{RHG}.

\begin{lem}\label{intHl}
For any $t\in G$, and any distinct $\alpha , \beta \in \Lambda $,
the intersection $H_\alpha \cap H_\beta ^t$ is finite. In
particular, if $a\in H_\alpha $ and $b\in H_\beta $ are elements of
infinite order, then $a$ and $b$ are not commensurable in $G$.
\end{lem}

\begin{proof}[Proof of Corollary \ref{frh}]
We assume that $G$ is non--elementary and $H_\lambda $ is proper for
any $\lambda \in \Lambda $. (Otherwise the corollary is obvious.) By
Lemma \ref{hypel} there is a hyperbolic element of infinite order
$g\in G$. By Lemma \ref{E(g)}, $G$ is hyperbolic relative to
$\Hl\cup \{ E_G(g)\} $. Thus without loss of generality we may
assume that $H_{\lambda _1}$ is infinite elementary for some
$\lambda _1\in \Lambda $. Note that $H_{\lambda _1}$ is proper as
$G$ is non--elementary. Applying the same arguments again, we may
assume that there is $\lambda _2\in \Lambda $, $\lambda _2\ne
\lambda _1$, such that $H_{\lambda _2} $ is also infinite and
elementary.

Let $S$ be a subset of $G$ and let $\mathcal F=\mathcal F(S)$ be the
set provided by Theorem \ref{CEP}. For each $\lambda \notin \{
\lambda _1, \lambda _2\} $, let $N_\lambda $ be a normal subgroup in
$H_\lambda $ such that $H_\lambda /N_\lambda $ is hyperbolic and the
natural homomorphism $H_\lambda \to H_\lambda/N_\lambda $ is
injective on $(H_\lambda \cap \mathcal F)\cup \{ 1\} $.  We also set
$N_{\lambda _1}=N_{\lambda _2}=\{ 1\} $. In particular, $N_\lambda
\cap \mathcal F=\emptyset $ for all $\lambda \in \Lambda $.

By Theorem \ref{CEP} the group $G(\N )$ is hyperbolic relative to a
collection of hyperbolic subgroups, i.e., it is hyperbolic itself as
observed in Corollary \ref{hypquot}. Moreover, the restriction of
the natural homomorphism $G\to G(\N )$ to $S$ is injective. To show
that $G(\N )$ is non--elementary, it suffices to note that $G(\N )$
contains at least two non--commensurable elements of infinite order.
Indeed any two elements of infinite order $h_1\in H_{\lambda
_1}/N_{\lambda _1}\cong H_{\lambda _1}$ and $h_2\in H_{\lambda
_1}/N_{\lambda _2}\cong H_{\lambda _2}$ are not commensurable
according to Lemma \ref{intHl}.
\end{proof}

\vspace{1cm}

\noindent Denis Osin\\ Department of Mathematics \\ The City College of New York\\
138th street and Convent Ave.\\ New York, NY 10031\\

\noindent {\rm E-mail address:} \it denis.osin@gmail.com

\end{document}